\documentclass[twoside, reqno, 10pt]{amsart}

\usepackage[left=2.5cm, right=2.5cm, top=3cm, bottom=2.5cm]{geometry}

\usepackage[colorlinks=true, pdfstartview=FitV, linkcolor=blue,citecolor=blue, urlcolor=blue]{hyperref}

\usepackage[usenames]{xcolor} %I put this in to get rid of the incompatible color warnings -ken

\definecolor{labelkey}{rgb}{0,0,1}
\definecolor{Red}{rgb}{0.7,0,0.1}
\definecolor{Green}{rgb}{0,0.7,0}

\usepackage{amsfonts, amssymb, amsmath, amsthm, mathrsfs, bbm, cjhebrew, gensymb, textcomp, mathtools, dsfont, calligra}

\usepackage[normalem]{ulem}%normalem gets rid of underlining of emph

\usepackage{todonotes}
\usepackage{esint}

\usepackage{commath, setspace, subcaption, parcolumns, multirow, multicol, accents, comment, marginnote, verbatim, empheq, enumerate, stackrel, enumitem, float, physics}
\usepackage[capitalize,nameinlink,noabbrev]{cleveref}

\usepackage{graphicx, graphics, epsfig, psfrag, tikz, tikz-cd,svg}

\numberwithin{equation}{section}

\usepackage{tikz}

\newcommand{\etalchar}[1]{$^{#1}$}
\providecommand{\bysame}{\leavevmode\hbox to3em{\hrulefill}\thinspace}
\providecommand{\MR}{\relax\ifhmode\unskip\space\fi MR }
\providecommand{\href}[2]{#2}

\newtheorem{Thm}{Theorem}[section]
\newtheorem{Lem}[Thm]{Lemma}

\newtheorem{Cor}[Thm]{Corollary}

\newtheorem{Rmk}[Thm]{Remark}

\newtheorem*{Thm*}{Theorem}

\newcommand{\Z}{\mathbb{Z}}

\newcommand{\C}{\mathbb{C}}
\newcommand{\T}{\mathbb{T}}

\def\AA {\mathcal{A}}

\newcommand{\PP}{\mathcal{P}}
\newcommand{\QQ}{\mathcal{Q}}
\newcommand{\FF}{\mathcal{F}}
\newcommand{\EE}{\mathcal{E}}
\newcommand{\ZZ}{\mathcal{Z}}

\newcommand{\DD}{\mathcal{D}}

\newcommand{\II}{\mathcal{I}}
\newcommand{\JJ}{\mathcal{J}}
\newcommand{\KK}{\mathcal{K}}
\newcommand{\LL}{\mathcal{L}}
\newcommand{\RR}{\mathcal{R}}

\def\SS {\mathcal{S}}
\newcommand{\tKK}{\tilde{\KK}}
\newcommand{\tLL}{\tilde{\LL}}

\newcommand{\tQQ}{\til{\mathcal{Q}}}

\newcommand{\PPd}{\dot{\mathcal{P}}}
\newcommand{\ZZd}{\dot{\mathcal{Z}}}
\newcommand{\EEd}{\dot{\mathcal{E}}}
\newcommand{\JJd}{\dot{\mathcal{J}}}

\newcommand{\EEdd}{\ddot{\mathcal{E}}}

\DeclareMathOperator{\esssup}{ess\ sup}

\newcommand{\til}[1]{{\tilde{#1}}}
\newcommand{\Sob}[2]{\lVert#1\rVert_{#2}}

\newcommand{\goesto}{\rightarrow}
\newcommand{\smod}{\setminus}

\newcommand{\al}{\alpha}
\newcommand{\be}{\beta}
\newcommand{\de}{\delta}
\newcommand{\De}{\Delta}
\newcommand{\gam}{\gamma}

\newcommand{\eps}{\epsilon}
\newcommand{\veps}{\varepsilon}

\newcommand{\s}{\sigma}
\newcommand{\lam}{\lambda}
\newcommand{\kap}{\kappa}

\newcommand{\w}{\omega}

\newcommand{\bdy}{\partial}
\newcommand{\lb}{\langle}
\newcommand{\rb}{\rangle}

\newcommand{\tu}{\tilde{u}}
\newcommand{\tnu}{\tilde{\nu}}
\newcommand{\tG}{\tilde{G}}

\newcommand{\tal}{\tilde{\al}}
\newcommand{\tc}{\tilde{c}}
\newcommand{\tU}{\tilde{U}}

\DeclareRobustCommand\nlab{k}
 
\begin{document}

\title[Viscosity Parameter Estimation for 2D NSE]{Convergence Analysis of a Viscosity Parameter Recovery Algorithm for the 2D Navier-Stokes Equations}

\author{Vincent R. Martinez}

\thanks{}

%\date{\today}

\maketitle

\begin{abstract}
In this paper, the convergence of an algorithm for recovering the unknown kinematic viscosity of a two-dimensional incompressible, viscous fluid is studied. The algorithm of interest is a recursive feedback control-based algorithm that leverages observations that are received continuously-in-time, then dynamically provides updated values of the viscosity at judicious moments. It is shown that in an idealized setup, convergence to the true value of the viscosity can indeed be achieved under a natural and practically verifiable non-degeneracy condition. This appears to be first such result of its kind for parameter estimation of nonlinear partial differential equations. Analysis for two parameter update rules is carried out: one which involves instantaneous evaluation in time and the other, averaging in time. The proofs of convergence for either rule exploits sensitivity-type bounds in higher-order Sobolev topologies, while the instantaneous version particularly requires delicate energy estimates involving the time-derivative of the sensitivity-type variable. Indeed, a crucial component in the analysis of the first update rule is the identification of a dissipative structure for the time-derivative of the sensitivity-type variable, which ultimately ensures a favorable dependence on the tuning parameter of the algorithm.
\end{abstract}

\vspace{1em}

{\noindent \small {\it {\bf Keywords: parameter estimation, system identification, data assimilation, feedback control, nudging, synchronization, Navier-Stokes equations, viscosity, convergence, sensitivity}
  }  
  }

\section{Introduction}\label{sect:intro}

The main problem in consideration is the determination of the value of an apriori unknown material parameter of a dynamical system via a time-series of observations made on a subset of its corresponding state variables. This problem is referred to as the problem \textit{parameter estimation} and constitutes a fundamental difficulty in the development of predictive models in the physical sciences, engineering, and social sciences. In the particular field of fluid dynamics, for instance, for groundwater, the transmissivity coefficient of the water, which captures its ability to move to across an aquifer, is a central parameter in these models that must be estimated in some way \cite{Yeh1986, EwingLin1991}. Significant efforts in the past several decades have been dedicated to developing techniques for inferring model parameters of nonlinear dynamical systems based on knowledge of the system itself, a time-series of observations on the state variables, and possibly the statistics of the noise that may corrupt the observations or the model (see, for example, \cite{Cialenco2018, AbarbanelCrevelingFarsianKostuk2009, VossTimmerKurths2004} and references therein). In spite of these developments, rigorous analytical results for parameter estimation algorithms applied to inferring parameters of nonlinear partial differential equations (PDEs) have appeared only relatively recently \cite{CialencoGlattHoltz2011, CarlsonHudsonLarios2020}. In these works, results were mostly limited to demonstrating consistency of the estimator, asymptotic normality, or performing a sensitivity-type analysis. The main purpose of this paper is to address the issue of convergence of one such algorithm as applied to the problem of estimating the kinematic viscosity of a two-dimensional (2D) incompressible fluid in an ideal setup.

The particular system of interest is the forced, two-dimensional (2D) Navier-Stokes equations (NSE) over the periodic box $\T^2=[0,2\pi]^2$. Given a sufficiently smooth, mean-free vector field $f$, a sufficiently smooth, mean-free initial velocity field $u_0$, both of which are $2\pi$-periodic in $x,y$, this system is given by
    \begin{align}\label{eq:nse}
        \bdy_tu+(u\cdotp\nabla)u=\nu\De u-\nabla p+f,\quad \nabla\cdotp u=0,
    \end{align}
where $\nu>0$ denotes the kinematic viscosity and $p$ denotes the scalar pressure field, which is determined by $u$ through the equation $\De p=\sum_{i,j}\bdy_i\bdy_j(u^iu^j)$. Upon passing to frequency-variables via the Fourier transform, one may equivalently view \eqref{eq:nse} as an infinite-dimensional dynamical system governing the evolution of each Fourier mode of the flow field $u$. In adopting this point of view, we will consider observations on the flow field $u$ to be given by a continuous time-series, $\{P_Nu(t)\}_{t\geq0}$, for some $N>0$, where $P_N$ denotes the low-pass filter that projects vector fields onto the subspace determined by their Fourier modes through wave-number $|k|\leq N$. We will refer to such observations as ``spectral observations." Note that our setting is highly idealized in that we assume a perfect model, perfect observations, with those observations made continuously-in-time and presented in a spectral form. This allows us to isolate the difficulty brought on by the infinite-dimensionality of \eqref{eq:nse} alone in establishing the desired convergence result. We refer the reader to \cref{rmk:real} for additional discussion on more physical situations to be studied in subsequent works.

In \cite{CarlsonHudsonLarios2020}, an algorithm was introduced for dynamically recovering the unknown kinematic viscosity of \eqref{eq:nse} given spectral observations of the flow-field collected continuously-in-time. The algorithm relies on a certain proxy, $\tu$, for the full flow-field that is generated by integrating a feedback-controlled system, in which the spectral observations are directly inserted, forward in time. Given
$\tnu>0$, $\mu>0$, and $N\geq1$, this system is given by
    \begin{align}\label{eq:nse:ng}
        \bdy_t\tu+(\tu\cdotp\nabla)\tu=\tnu\De\tu-\nabla q+f-\mu P_N(\tu-u),\quad \nabla\cdotp\tu=0,
    \end{align}
where $q$ denotes the corresponding scalar pressure field, $\mu$ represents a tuning parameter, and $\tnu$ denotes some representative value, i.e., a ``guess", for the true kinematic viscosity $\nu$. Denoting the difference in flow fields by $w=\tu-u$, it was proposed in \cite{CarlsonHudsonLarios2020} that a new representative value $\nu$ be generated according to the following formula:
    \begin{align}\label{def:new:rep}
        \nu_{new}=\tnu+\mu\frac{\Sob{P_Nw(t)}{L^2}^2}{\lb P_N(-\De)\tu(t),P_Nw(t)\rb_{L^2}}\bigg{|}_{t=t_*},
    \end{align}
where evaluation occurs at some judiciously chosen time $t_*>0$. Since $P_Nw$ depends only on the given data, $P_Nu$, and the computed proxy, $\tu$, the right-hand side can indeed be evaluated. Upon generating, $\nu_{new}$, the system \eqref{eq:nse:ng} with $\tnu=\nu_{new}$ is then integrated forward-in-time once again in order to generate an updated representative value via \eqref{def:new:rep}, and so on. Several numerical simulations were carried out in \cite{CarlsonHudsonLarios2020} that tested the efficacy of this algorithm in a regime for turbulent flows. For sufficiently large values of the tuning parameter $\mu>0$ and number of observed modes $N>0$, convergence was reliably observed. In this paper, we provide analytical confirmation of these observations provided that \eqref{def:new:rep} is well-defined at each instance that it is evaluated and that a certain non-degeneracy condition holds. Indeed, the problem of inferring the true value of $\nu$ given knowledge of $P_Nu(t)$, for some $N>0$ and for all $t\geq0$ may be viewed as an ill-posed inverse problem. Morally speaking, the non-degeneracy condition that is identified here constrains one to a ``well-posed" regime.

The inspiration for this algorithm derives from a seminal work of Azouani, Olson, and Titi \cite{AzouaniOlsonTiti2014}. There, a beautifully simple algorithm for data assimilation for systems of partial differential equations based on suitably interpolating observations into the phase space of the system, then inserting them into the model as a feedback control as in \eqref{eq:nse:ng} above, was proposed and subsequently studied for the 2D NSE as a paradigmatic example. It was shown that when $\tnu=\nu$, then for a sufficiently large density of observations, one may identify a range of values for the tuning parameter, $\mu$, that guarantees synchronization of the approximating signal, $\tu$, to the reference signal, $u$, in the time-infinite limit. Since then, this synchronization-based approach for data assimilation of PDEs has been studied in a plethora of works that provide rigorous confirmation and justification, in the form of mathematical proof, for various practices in data assimilation such as reduced-component assimilation \cite{FarhatJollyTiti2015, FarhatLunasinTiti2016a}, assimilation of only temperature measurements \cite{FarhatLunasinTiti2016b, FarhatLunasinTiti2016c, FarhatGlattHoltzMartinezMcQuarrieWhitehead2020},  and assimilation of only boundary data \cite{JollyMartinezTiti2017}, to name only a few. This approach has also been exploited to obtain uniform-in-time error bounds for a post-processing Galerkin method \cite{MondainiTiti2018}, developed for the downscaling of statistical solutions of the 2D NSE \cite{BiswasFoiasMondainiTiti2018}, extended to accommodate local observations \cite{BiswasBradshawJolly2020}, as well as found to exhibit strong convergence properties \cite{BiswasMartinez2017, BiswasBrownMartinez2021}. Recently, studies in how well the approach reconstructs turbulent flows have also been carried out \cite{BuzzicottiClarkDiLeoni2020, ClarkDiLeoniMazzinoBiferale2020}. For a more comprehensive list of recent works on this particular approach to data assimilation, we refer the reader to \cite{CarlsonHudsonLariosMartinezNgWhitehead2021}.

In the setting where the model or observations may be corrupted by noise, a classical Kalman-filter based approach for parameter estimation can be applied by simply extending the state space of the system to include the parameters as additional states. This was originally introduced in \cite{Mayne1963} in the setting of linear systems (see also \cite{FreemanHassanMorton1986}). In an attempt to develop a connection, consider the (discrete) system dynamics is determined by $u_{k+1}=Fu_k+\xi_k$, for each time step $k\geq0$, where $F$ denotes the \textit{unknown} matrix of constants, and $\xi_k$ represents random model noise, i.e., observations are measured perfectly. Assuming gaussianity of $\xi$, the ``optimal parameter estimate" \cite{Mayne1963} is then given by
    \begin{align}\notag
        \til{F}_{k+1}=\til{F}_k+\kap_k(u_k-\til{u}_k),
    \end{align}
where $\kap_k$ is the so-called Kalman gain matrix and $\til{u}_k$ denotes the ``optimal state estimate." Although $\kap_k$ depends on the  conditional error variance associated to the optimal estimator, $\til{F}$, one may nevertheless na\"ively derive some intuition for the update formula given by \eqref{def:new:rep}. Traditionally, $\iota_k:=u_k-\til{u}_k$ is referred to as the \textit{innovation}. Thus, in the case of \eqref{def:new:rep}, one has the following correspondence  
    \[
        \kap_k\sim \frac{\mu}{\lb P_N(-\De)\tu(t_k),P_N(u(t_k)-\tu(t_k))\rb_{L^2}},\quad \iota_k\sim \Sob{P_N(u(t_k)-\tu(t_k))}{L^2}^2.
    \]
It would be interesting to make this apparent connection more rigorous, but do not pursue this point further here. We refer the reader to  \cite{LawStuartZygalakisBook} for a presentation of other filter-based techniques for data assimilation and the review article \cite{SchoukensLjung2019} for a more recent exposition on approaches to the problem of nonlinear system identification.

The proof of convergence of \eqref{def:new:rep} presented here was inspired by the recent work of the author and collaborators \cite{CarlsonHudsonLariosMartinezNgWhitehead2021}. There, convergence of the algorithm was analogously studied for the Lorenz equations. The crucial ingredient was in identifying a Lyapunov-type structure for the time-derivative of the system. Remarkably, we show that this approach can also be adapted to the case of the 2D NSE, although the technicalities involved are far greater. Indeed, the Lyapunov-type structure identified here is more sophisticated and requires a more nuanced analysis to properly exploit. In a word, the analysis of time-derivatives naturally introduce higher-order derivatives in the estimates. This, in turn, requires access to higher-order bounds for both the approximating signal, $\tu$, and the ``pseudo-sensitivity" variable, $w=\tu-u$. The Lyapunov-type structure must ultimately be exploited in a rather delicate way to close the estimates. In the setting of finite-dimensional dynamical systems, we also refer the reader to the works of \cite{AbarbanelCrevelingJeanne2008, CrevelingJeanneAbarbanel2008,  QuinnBryantCrevelingKleinAbarbanel2009} for various studies on similar synchronization-based approach to parameter estimation. We refer the reader as well to \cite{CotterDashtiRobinsonStuart2009} for a Bayesian approach to inverse problems in fluid mechanics. In the particular setting of finite-dimensional systems, we also point out that the convergence analysis performed in \cite{CarlsonHudsonLariosMartinezNgWhitehead2021} appears to be the first such result of its kind. We refer the reader to \cref{rmk:lorenz} for additional comments on how the case treated there compares to the one treated here. In the course of writing this article, the author was informed of a forthcoming paper that also studies the parameter estimation problem for the 2D NSE, but through variational techniques, where convergence-type results are also obtained \cite{BiswasHudson2021}.

We will also present a convergence analysis for a ``time-averaged" version of the viscosity update rule \eqref{def:new:rep} given by
    \begin{align}\label{def:mod:rep}
        \bar{\nu}_{new}=\tnu+\frac{\mu\frac{1}{t-t'}\int_{t'}^t\Sob{P_Nw(s)}{L^2}^2ds}{\frac{1}{t-t'}\int_{t'}^t\lb P_N(-\De)\tu(s),P_Nw(s)\rb_{L^2} ds}\bigg{|}_{t=t_*,t'=t_*'},\quad t_*>t_*',
    \end{align}
where, again, $t_*, t_*'$ is judiciously chosen. Ultimately, we establish an analogous convergence result for the algorithm induced by \eqref{def:mod:rep}. In comparison with the proof for \eqref{def:new:rep}, the proof for \eqref{def:mod:rep} will be significantly easier due to the fact that control of time-derivatives of $\Sob{w_N(t)}{L^2}^2$ is not required. Given the apparent flexibility available by for inferring viscosity that is suggested by \eqref{def:mod:rep} and \eqref{def:new:rep}, one may be naturally inclined to ask whether this nudging-based approach can be used to infer parameters other than the viscosity in \eqref{eq:nse} such as the external forcing $f$. This interesting and very important issue will be explored in future works.

The organization of this paper is as follows. In \cref{sect:background}, we develop the mathematical setting and notation in order to provide precise statements of the results.  The main results are then stated \cref{sect:statements}. In \cref{sect:proof}, the proofs of the main convergence results are supplied. A considerable effort of this paper is dedicated to establishing the time-derivative estimates required to prove \cref{thm:converge}. The paper concludes with \cref{sect:app:apriori} and \cref{sect:app:sensitivity}, where higher-order estimates for $\tu$ and $w$ are respectively established.

\section{Mathematical Setting and Notation}\label{sect:background}
Let $\mathcal{B}(\T^2)$ denote the set of all Borel measureable functions on $\T^2$ that are $2\pi$-periodic a.e. in each direction. Then for $1\leq p\leq\infty$, the Lebesgue spaces are defined by 
    \begin{align}\label{def:Lp}
        L^p(\T^2):=\{\phi\in\mathcal{B}(\T^2):\Sob{\phi}{L^p}<\infty\},
    \end{align}
where
    \begin{align}\label{def:Lp:norm}
    \Sob{\phi}{L^p}^p=\begin{cases}\int_{\T^2}|\phi(x)|^pdx\\
    \esssup_{x\in\T^2}|\phi(x)|.
    \end{cases}
    \end{align}
For $s\geq0$, we define the (homogeneous) Sobolev spaces by
    \begin{align}\label{def:Sob:sp:hom}
        H^s(\T^2):=\{\phi\in L^2(\T^2):\Sob{\phi}{H^s}<\infty,\ \hat{\phi}(\mathbf{0})=0\},
    \end{align}
where $\mathbf{0}=(0,0)\in\Z^2$, and
    \begin{align}\label{def:Sob:hom}
        \Sob{\phi}{H^s}^2=\sum_{k\in\Z^2}|k|^{2s}|\hat{\phi}(k)|^2,
    \end{align}
where $\hat{\phi}(k)$ denotes the Fourier coefficient of the Fourier series of $\phi$ corresponding to wave-number $k$. By the definition of the Fourier transform, we see that
    \begin{align}\label{eq:mean:zero}
        \int_{\T^2}\phi(x)dx=\hat{\phi}({\mathbf{0}})=0,
    \end{align}
whenever $\phi\in H^s$. In particular $H^s(\T^2)$ consists of zero-mean functions. By Plancherel's theorem, we further see that 
    \begin{align}\label{eq:Sob:char}
        \Sob{\nabla \phi}{L^2}=\Sob{\phi}{H^1}^2,\quad \Sob{\phi}{H^s}=\Sob{(-\De)^{s/2}\phi}{L^2},\quad s\geq0.
    \end{align}
In this setting, one has a Poincar\'e inequality, which implies
    \begin{align}\label{eq:poincare}
        \Sob{\phi}{L^2}\leq \Sob{\phi}{H^s},
    \end{align}
for any $s\geq0$. In particular, it immediately follows that
    \begin{align}\label{eq:hom:inhom:equiv}
        \Sob{\phi}{H^s}^2\leq \sum_{k\in\Z^2}(1+|k|^2)^s|\hat{\phi}(k)|^2\leq 2^s\Sob{\phi}{H^s}^2.
    \end{align}
The intermediate quantity is usually referred to as the (square of the) inhomogeneous Sobolev norm of $\phi$. The corresponding 2D vector field analogs of the Lebesgue and Sobolev spaces are then simply $(L^p(\T^2))^2$ and $(H^s(\T^2))^2$. For convenience, we adopt an abuse of the notation and do not distinguish between the Lebesgue and Sobolev spaces for scalar functions or vector fields. We will also adopt the following notation for $L^2$--norms and $H^1$--norms:
    \begin{align}\label{def:foias:temam}
        |\phi|=\Sob{\phi}{L^2},\quad \lVert\phi\rVert=\Sob{\phi}{H^1}.
    \end{align}

We then define the corresponding Lebesgue and Sobolev spaces of divergence-free vector fields by
    \begin{align}
        L^2_\s(\T^2)&=\{v\in L^2(\T^2): |v|<\infty,\ k\cdotp\hat{v}(k)=0,\ k\neq0\},\label{def:L2:sol}
    \end{align}
and
    \begin{align}
         H^s_\s(\T^2)&=L^2_\s(\T^2)\cap H^s(\T^2).\label{def:Hs:sol}
    \end{align}

It is customary to handle the pressure term in \eqref{eq:nse} by projecting the system onto the subspace of divergence-free vector fields in $L^2$. This is done via the Leray projection, $P_\s$, which is defined by
    \begin{align}\label{def:leray}
        \left(\widehat{P_\s v}\right)^i(k)=\left(1-\frac{k_ik_j}{|k|^2}\right)\hat{v}^j(k),\quad k=(k_1,k_2)\in\Z^2\smod\{{\mathbf{0}}\},\quad i=1,2.
    \end{align}
We see that in this setting of periodic boundary conditions, $P_\s$ commutes with spatial derivatives. Hence
    \begin{align}\notag
        |A^{s/2}\phi|\leq\Sob{\phi}{H^s}\leq 2^{s/2}|A^{s/2}\phi|.
    \end{align}
Upon setting $g=P_\s f$, then formally applying $P_\s$ to \eqref{eq:nse} and \eqref{eq:nse:ng}, we have
    \begin{align}
        \frac{du}{dt}+\nu Au+B(u,u)&=g,\label{eq:nse:ff}\\
        \frac{d\tu}{dt}+\tnu A\tu+B(\tu,\tu)&=g-\mu P_N(\tu-u),\label{eq:nse:ng:ff}
    \end{align}
where $A=-P_\s\De$ denotes the Stokes operator and 
    \begin{align}\label{def:B}
        B(u,u):=P_\s[(u\cdotp\nabla)u].
    \end{align}
Then \eqref{eq:nse:ff}, \eqref{eq:nse:ng:ff} is referred to as the functional form of \eqref{eq:nse}, \eqref{eq:nse:ng}, respectively. We will make use of the following important orthogonality property of the bilinear form:
    \begin{align}\label{eq:ortho:B}
        \lb B(u,v),v\rb=0=\lb B(u,u),Au\rb,
    \end{align}
where the first identity holds for $u,v\in H^1_\s(\T^2)$ and the second identity holds for $u\in H^2_\s(\T^2)$.

Given $g\in L_\s^2(\T^2)$, we define the Grashof number by
    \begin{align}\label{def:Grashof}
        G:=\frac{|g|}{\nu^2}.
    \end{align}
One then has the following classical results regarding the existence theory for \eqref{eq:nse} (see \cite{ConstantinFoias1988, TemamBook1997, TemamBook2001, FoiasManleyRosaTemamBook2001}).

\begin{Thm}\label{thm:nse:wellposed}
Let $g\in H^{k-1}_\s(\T^2)$, for some $k\geq1$. Then for all $u_0\in H^k_\s(\T^2)$, there exists a unique vector field $u\in C([0,T];H_\s^k(\T^2))\cap L^2(0,T;H^{k+1}(\T^2))$ satisfying \eqref{eq:nse:ff}, for all $T>0$, and the estimate
    \begin{align}\label{est:nse:H1}
        \lVert u(t)\rVert^2\leq e^{-\nu (t-t_0)}\lVert u(t_0)\rVert^2+\nu^2G^2(1-e^{-\nu(t-t_0)}),
    \end{align}
for all $t\geq t_0\geq0$. Moreover, there exists an open neighborhood $U\subset\C$ such that $U\supset (0,\infty)$ and $u:U\goesto H^2_\s(\T^2)$ is analytic.
\end{Thm}

From \cref{thm:nse:wellposed}, we see that the operator $S_g(t):H^k_\s(\T^2)\goesto H^k_\s(\T^2)$, $S_g(t)u_0=u(t;u_0,g)$, where $u(t;u_0,g)$ denotes the unique solution of \eqref{eq:nse:ff} corresponding to $g\in L^2_\s(\T^2)$ and $u_0\in H^k_\s(\T^2)$, evaluated at time $t\geq0$, defines a continuous semigroup (see \cite{TemamBook1997}). From \eqref{est:nse:H1}, we see that for $t_0=0$ and $t>0$ sufficiently large, one has
    \begin{align}\label{def:rad:H1}
        \lVert u(t)\rVert\leq \sqrt{2}\nu G=:R_1.
    \end{align}
This defines an absorbing ball for the semigroup, $\{S(t)\}_{t\geq0}$, of \eqref{eq:nse:ff} with respect to the $H^1$--topology. We denote the ball of radius $R_1$, centered at the origin in $H^1(\T^2)$, by $B_1(R_1)$. Then one has $S(t)u_0\in B_1(R_1)$, for all $u_0\in B_1(R_1)$.

When $g\in H^{k-1}_\s(\T^2)$, for some $k\geq2$, it was shown in \cite{DascaliucFoiasJolly2005}, for the case $k=2$, and  \cite{BiswasBrownMartinez2021}, for the general case of $k\geq3$, that \eqref{eq:nse:ff} has an absorbing ball, $B_k(R_k)$, with respect to the $H^k$--topology whose radius is given by
    \begin{align}\label{def:rad:Hk}
        |A^{k/2}u(t)|\leq c_k\nu(\s_{k-1}^{1/k}+G)^{k-1}G=:R_k,
    \end{align}
for some constant $c_k>0$, depending only on $k$, where $\s_\ell$ denotes higher order ``shape factors" of the force defined by
    \begin{align}\label{def:shape}
        \s_\ell:=\frac{|A^{\ell/2}g|}{|g|}.
    \end{align}
Observe that $\s_\ell\geq1$, for all $\ell\geq0$, by \eqref{eq:poincare}.

The basic existence theory for \eqref{eq:nse:ng:ff} was initially developed in \cite{AzouaniOlsonTiti2014} in a rather general setting of interpolant observable quantities that includes the case of spectral observables, which are the particular form of observables considered in this article.  More recently, the existence theory has been extended to accommodate a more general class of interpolant observable quantities in order to study the nature of synchronization is spaces of arbitrarily large Sobolev regularity in \cite{BiswasBrownMartinez2021}. In the particular case of spectral observations though, the existence theory and topology of synchronization was actually extended to the setting of analytic Gevrey classes in \cite{BiswasMartinez2017}. Also, the analogous statement regarding time analyticity for solutions to the feedback control system was established in \cite{IbdahMondainiTiti2019}. In all of the above works, only the case $\nu=\tnu$ was treated. The case where $\nu\neq\tnu$ was initially studied in \cite{CarlsonHudsonLarios2020}. %We state the following existence result, which follows from trivial modifications of the original proofs found in \cite{AzouaniOlsonTiti2014}, \cite{IbdahMondainiTiti2019}.

\begin{Thm}
Let $g\in H_\s^{k-1}(\T^2)$, for some $k\geq1$. Suppose $u_0\in H_\s^k(\T^2)$ and that $\nu,\tnu>0$. Let $u(\cdotp;u_0,\nu)$ denote the unique, global-in-time solution of \eqref{eq:nse:ff} corresponding to initial data $u_0$ and viscosity $\nu$. There exists a constant $\tc_0\geq1$ such that if $\mu>0$ and $N\geq1$ satisfy
    \begin{align}\label{cond:mu:N:tnu}
        \mu\leq \tc_0N^2\tnu,
    \end{align}
then for all $\tu_0\in H_\s^k(\T^2)$, there is a unique  $\tu \in C([0,T];H_\s^k(\T^2))\cap L^2(0,T;H^{k+1}(\T^2))$ satisfying \eqref{eq:nse:ng:ff}, for all $T>0$. Moreover, there exists an open neighborhood $\tU\subset\C$ such that $\tU\supset(0,\infty)$ and $\tu:\tU\goesto H^2_\s(\T^2)$ is analytic.
\end{Thm}

In order to prove convergence of the parameter update algorithm described in \cref{sect:intro}, we will ultimately make use of bounds analogous to \eqref{def:rad:H1}, \eqref{def:rad:Hk} for the solutions to \eqref{eq:nse:ng:ff}. Since such bounds do not seem to be available in the existing literature, we establish them here. The statement of this result will be postponed for \cref{sect:statements}, while its proof will be supplied in \cref{sect:app:apriori}.

\section{Statements of Main Theorems}\label{sect:statements}

Our first main result (\cref{thm:ng:bounds}) asserts bounds for solutions to \eqref{eq:nse:ng:ff} that carefully track the dependence on the tuning parameter $\mu$. Such bounds have been obtained in \cite{AzouaniOlsonTiti2014} in the case when $\nu=\tnu$ and \cite{CarlsonHudsonLarios2020} in the case when $\nu\neq\tnu$. However, these bounds were only obtained in $L^2(\T^2)$ and $H^1(\T^2)$, whereas we will also need them in $H^2(\T^2)$. Most importantly, our convergence analysis is dynamical and will require $\tnu$ to change values. We must, therefore, carefully track the dependence on $\tnu$ relative to $\nu$. Lastly, since our convergence result is asymptotic in nature, it will be convenient to frame our results in a sufficiently evolved regime for \eqref{eq:nse:ff}. In particular, we will typically assume that $u_0\in B_k(R_k)$, for some $k\geq1$, where $R_k$ is given by \eqref{def:rad:Hk}, while $\tu_0$ will belong to the same ball of a possibly inflated radius. Due to these various adjustments, we find it necessary to re-prove these estimates to suit our purposes. Lastly, we point out that although our analysis may not be the sharpest possible in the setting of periodic boundary conditions treated here, we give particular emphasis to the tracking of the dependence on all relevant system parameters and the manner in which they scale in the estimates. %We refer the reader to \cite{AzouaniOlsonTiti2014}, where the best-known bounds on $\mu, N$ to guarantee synchronization of the nudging-based scheme can be found and  We also refer the reader to \cref{rmk:algorithm} for relevant related remarks.

Before we state these bounds, let us first introduce a modified Grashof number, $\tG$, defined by
    \begin{align}\label{def:tG}
        \tG:=\left(\left(\frac{\nu}{\tnu}\right)\left(\frac{\nu}{\mu}\right) +1\right)^{1/2}G.
    \end{align}
Observe that $G\leq\tG$. 

\begin{Thm}\label{thm:ng:bounds}
Let $g\in H_\s^{k-1}(\T^2)$, for some $k\geq1$, $u_0\in B_\ell(R_\ell)$, and $\tu_0\in B_\ell(\al R_\ell)$, for all $\ell=1,\dots, k$, for some $\al\geq1$, where $R_\ell$ is given by \eqref{def:rad:Hk}. Suppose that $\mu, N, \tnu>0$ satisfy \eqref{cond:mu:N:tnu}. Let $u(\cdotp;u_0,\nu), \tu(\cdotp;\tu_0,\tnu)$ denote the corresponding unique global-in-time solutions to \eqref{eq:nse:ff}, \eqref{eq:nse:ng:ff} with initial values $u_0, \tu_0$ and viscosities $\nu,\tnu$, respectively. There exists a constant $c_0\geq1$, such that if $\mu$ additionally satisfies
    \begin{align}\label{cond:mu:N:nu}
        \mu\leq c_0N^2\nu,
    \end{align}
then the following statements hold: If $k=1$, then there exists a constant $\tal_1\geq1$, depending on $\al$, such that
    \begin{align}\label{def:rad:H1:ng}
        \sup_{t\geq0}\lVert\tu(t)\rVert^2\leq \tal_1^2\nu^2\tG^2.
    \end{align}
If $k=2$, then there exist constants $\al_2,\tal_2\geq1$, depending on $\al$, such that if $\mu$ additionally satisfies
    \begin{align}\label{cond:mu:ng:H2}
        \mu\geq \nu\al_2^2\left(\frac{\nu}{\tnu}\right)\tG^2,
    \end{align}
then
    \begin{align}\label{def:rad:H2:ng}
    \sup_{t\geq0}|A\tu(t)|\leq \nu\tal_2\s_{1}^{1/2}\left(\s_{1}^{1/2}+G\right)\tG.
    \end{align}
If $k\geq3$, then there exist constant $\al_k,\tal_k\geq1$, depending on $\al$, such that if $\mu$ additionally satisfies
\begingroup\makeatletter \def\maketag@@@#1{\hbox{\m@th\normalfont#1$_{\nlab}$}} 
    \begin{align}\label{cond:mu:ng:Hk}
        \mu\geq \nu\al_k^2
        \left[\left(\tal_1^2+\tal_2^2\right)\left(\s_1^{1/2}+\tG\right)\tG+\tal_2^{2/k}\left(\frac{\nu}{\tnu}\right)^{1-2/k}\left(\frac{\s_{k-1}^{1/k}+G}{G}\right)^{2/k}\frac{\s_1^{1/2}+G}{\s_{k-1}^{1/k}+G}\right],
    \end{align}
\endgroup
then
    \begin{align}\label{def:rad:Hk:ng}
    \sup_{t\geq0}|A^{k/2}\tu(t)|\leq \tal_k\s_{k-1}^{1/k}\left(\s_{k-1}^{1/k}+G\right)^{k-1}\tG.
    \end{align}
\end{Thm}

\begin{comment}
\begin{Thm}\label{thm:ng:bounds}
Let $g\in H_\s^{k-1}(\T^2)$, for $k=1,2$, and $u_0\in B_k(R_k)$, where $R_k$ is given by \eqref{def:rad:Hk}. Let $u$ denote the corresponding unique global-in-time solution to \eqref{eq:nse:ff}. Suppose that $\mu, N, \tnu$ satisfy \eqref{cond:mu:N:tnu} and
    \begin{align}\label{cond:mu:N:tnu}
        \mu\leq cN^2\tnu.
    \end{align}
Suppose also that $\tu_0\in B_k(\al R_k)$, for some $\al\geq1$. Let $\tu(\cdotp;\tu_0,g,u)$ denote the unique global-in-time solution to \eqref{eq:nse:ng:ff}. If $k=1$, then there exists a constant $\tal_1\geq1$, depending on $\al$, such that
    \begin{align}\label{def:rad:H1:ng}
        \lVert\tu(t)\rVert^2\leq \tal_1^2\nu^2\tG^2,
    \end{align}
for all $t\geq0$. If $k\geq2$, then there exist constants $\al_k,\tal_k\geq1$, depending on $\al$, such that if $\mu$ additionally satisfies
    \begin{align}\label{cond:mu:ng:H2}
        \mu\geq 100\nu\begin{cases}
        \al_2^2\left(\frac{\nu}{\tnu}\right)\tG^2,&k=2\\
        C\left[\left(\tal_1^2+\tal_2^2\right)\left(\s_1^{1/2}+\tG\right)\tG+\tal_2^{2/k}\left(\frac{\nu}{\tnu}\right)^{1-2/k}\left(\frac{\s_{k-1}^{1/k}+G}{G}\right)^{2/k}\frac{\s_1^{1/2}+G}{\s_{k-1}^{1/k}+G}\right],&k\geq3,
        \end{cases}
    \end{align}
then
    \begin{align}\label{def:rad:Hk:ng}
    |A^{k/2}\tu(t)|\leq \tal_k\s_{k-1}^{1/k}\left(\s_{k-1}^{1/k}+G\right)^{k-1}G,
    \end{align}
for all $t\geq0$.
\end{Thm}
\end{comment}

The proof of \cref{thm:ng:bounds} will be provided in 
\cref{sect:app:apriori}. 

\begin{Rmk}\label{rmk:equiv}
If $\nu,\tnu$ satisfy
    \begin{align}\label{cond:nu:tnu:small}
        \frac{|\tnu-\nu|}{\nu}<\eps,
    \end{align}
for some $\eps\in(0,1)$, then
    \begin{align}\label{eq:nu:tnu:equiv}
        1-\eps\leq \frac{\tnu}{\nu}\leq 1+\eps.
    \end{align}
Also, if additionally, $\mu,\nu$ satisfy $\mu\geq\nu$, then 
    \[
            G\leq\tG\leq\left(\frac{2-\eps}{1-\eps}\right)^{1/2}G.
    \]
In particular, under these additional assumptions on $\mu,\nu,\tnu$, any dependence on $\tnu$ appearing in \eqref{cond:mu:N:tnu}, \eqref{def:rad:H1:ng}, \eqref{cond:mu:ng:H2}, \eqref{def:rad:H2:ng}, $\eqref{cond:mu:ng:Hk}_k$, \eqref{def:rad:Hk:ng} can be replaced by $\eps, \nu$, and all instances of $\tG$ can be replaced by $G$ up to a multiplicative constant depending only on $\eps$.
\end{Rmk}

The next result asserts bounds on the error between solutions of \eqref{eq:nse:ng:ff} and \eqref{eq:nse:ff} in terms of the difference between their respective viscosities, $\De\nu=\tnu-\nu$. Such bounds were obtained in \cite{CarlsonHudsonLarios2020}, but only in $L^2$ and $H^1$. In order to ultimately study the convergence of the algorithm induced by the formula \eqref{def:new:rep}, however, we must also have access to such bounds in $H^2$. It is, moreover, crucial to track the dependence of these bounds on the tuning parameter, $\mu$. We render this dependence in an explicit manner. More importantly, we require such estimate at the level of $H^2$.

\begin{Thm}\label{thm:sensitivity}
Let $\nu,\tnu,\tau_0>0$ and $g\in H_\s^2(\T^2)$. Suppose $u(\tau_0)\in B_\ell(R_\ell)$ and $\tu(\tau_0)\in B_\ell(\al R_\ell)$, for $\ell=1,2,3$, for some $\al\geq1$, where $R_3$ is given by \eqref{def:rad:Hk}. Suppose that $\mu, N,\nu,\tnu$ satisfy \eqref{cond:mu:N:tnu}, \eqref{cond:mu:N:nu}, \eqref{cond:mu:ng:H2}, and $\eqref{cond:mu:ng:Hk}_3$. Let $u(\cdotp;u(\tau_0),\nu)$, $\tu(\cdotp;\tu(\tau_0), \tnu)$ denote the unique global-in-time solutions to \eqref{eq:nse:ff}, \eqref{eq:nse:ng:ff} over $[\tau_0,\infty)$ corresponding to initial data $u(\tau_0), \tu(\tau_0)$, respectively. There exists a constants $c_1, c_2, c_3\geq1$, independent of $\al$, such that if $\mu$ additionally satisfies
    \begin{align}
        \mu&\geq c_1\nu\left[(\s_1^{1/2}+G)^2+(\s_2^{1/3}+G)^4\right]^{1/4}(\s_2^{1/3}+G)G,\label{cond:mu:sensitivity1}\\
        \mu &\geq c_2\nu\left(\frac{|\De\nu|}{\nu}\right)\tal_2^2\left(\s_1^{1/2}+\tG\right)^2\tG^2\label{cond:mu:sensitivity2},
    \end{align}
where $\tal_2$ is the constant from \cref{thm:ng:bounds}, then there exists $\tau>\tau_0$ such that
        \begin{align}\notag
    |A(\tu(t)-u(t))|^2\leq e^{-\mu(t-\tau)}|A(\tu(\tau)-u(\tau))|^2 +\nu^2\left(\frac{\nu}{\mu}\right)\left(\frac{|\De\nu|}{\nu}\right)^2K_2^2,
    \end{align}
holds for all $t\geq \tau$, where
    \begin{align}\label{def:K2}
        K_2^2:=c_3\tal_2^2\s_2^{2/3}(\s_2^{1/3}+G)^4\tG^2,
    \end{align}
and $\tal_2$ is the same constant appearing in \cref{thm:ng:bounds}.
\end{Thm}

\begin{Rmk}\label{rmk:higher:order}
It is not difficult to carry out the proofs of \cref{thm:ng:bounds} and \cref{thm:sensitivity} to arbitrarily high orders of Sobolev regularity. We refer the reader to \cite{FoiasGuillopeTemam1981, BiswasBrownMartinez2021} to carry out the argument to higher-order. In fact, since the observations are given as Fourier modes, one can establish the sensitivity-type bounds captured by \cref{thm:sensitivity} to analytic Gevrey classes provided that the external force $g$ belongs to such a class. We refer the reader to \cite{BiswasMartinez2017} for the relevant details.
\end{Rmk}

Our third result is a theorem that identifies conditions under which the formula \eqref{def:new:rep} converges to the true viscosity value. In order to make a precise statement, we must develop some notation: Given an integer $M\geq1$, suppose that a sequence of viscosities, $\{\nu_m\}_{m=0}^{M}$, a sequence of nudging parameters, $\{\mu_m\}_{m=1}^{M+1}$, and a sequence of times $0=t_0<\dots<t_{M}<\infty$ are given. Let $I_m:=[t_m,t_{m+1})$, for $0\leq m\leq M-1$. Given $\nu_0>0$, $g\in H^{k-1}_\s(\T^2)$, for some $k\geq1$, and $u_0,\tu_0\in H^k_\s(\T^2)$, we consider
    \begin{align}\label{eq:nse:ng:I0}
        \begin{split}
        \frac{du}{dt}+\nu Au+B(u,u)&=g,\quad u(0)=u_0,\\
        \frac{d\tu}{dt}+\nu_0 A\tu+B(\tu,\tu)&=g-\mu_1 P_N(\tu-u),\quad \tu(0)=\tu_0,
        \end{split}
    \end{align}
when $t\in I_0$. For $t\in I_m$, where $1\leq m\leq M$, we consider
    \begin{align}\label{eq:nse:ng:Im}
        \begin{split}
        \frac{du}{dt}+\nu Au+B(u,u)&=g,\quad u(t_m)=u(t_{m}^-;\nu),\\
        \frac{d\tu}{dt}+\nu_m A\tu+B(\tu,\tu)&=g-\mu_{m+1} P_N(\tu-u),\quad \tu(t_m)=\tu(t_{m}^-;\nu_{m-1}),
        \end{split}
    \end{align}
where $v(\tau^-)=\lim_{t\goesto\tau^-}v(t)$. In this framework, for $M\geq1$ and $1\leq m<M+1$, we define the following quantity:
    \begin{align}\label{def:eps}
        \veps_m:=|\lb AP_N\tu(t_m^-), P_N(\tu(t_m^-)-u(t_m^-))\rb|\nu_0^{-2}.
    \end{align}
Let us also suppose that each $\nu_m$ is related as follows:
    \begin{align}\label{def:nu:sequence}
      \nu_{m+1}=\nu_m+\mu\frac{\Sob{P_Nw(t)}{L^2}^2}{\lb P_N(-\De)\tu(t),P_Nw(t)\rb_{L^2}}\bigg{|}_{t=t_{m+1}},
    \end{align}
where $t_{m+1}>t_m$, for all $m=0,\dots, M-1$. 

\begin{comment}
\begin{Thm}\label{thm:converge}
Assume the hypotheses of \cref{thm:ng:bounds}.
Given $\nu_0\geq\nu$ and $M\geq1$, there exists a tuning parameter $\mu_0>0$, and a strictly increasing divergent sequence of update times $\{t_m\}_{m=0}^M$, such that $t_0=0$ and if
            \begin{align}\label{cond:non:degen}
              \veps:=\inf_{1\leq m<M+1}\veps_m>0
             \end{align}
for some $\veps>0$, then there exists $\be\in(0,1)$ such that for all $\mu\geq \mu_0$, there exists an observational threshold $N_0\geq1$ such that whenever $N\geq N_0$, one has
    \begin{align}\label{eq:converge}
        |\nu_{m+1}-\nu|\leq\be|\nu_m-\nu|,
    \end{align}
for all $m=0,\dots,M-1$, where $\nu_m$ is given by \eqref{def:new:rep}. Moreover, $\{\nu_m\}_{m=0}^M$ is a positive sequence of numbers.
\end{Thm}
\end{comment}

\begin{Thm}\label{thm:converge}
Given $\nu_0\geq\nu$ and $M\geq1$, there exist constants $\gam,\gam_1,\dots, \gam_M>0$, where $\gam$ depends only on the initial relative error, $\de_0:=(\nu_0-\nu)\nu^{-1}$, and the Grashof-type number $|g|\nu_0^{-2}$, and where $\gam_m$ depends additionally on $\veps_m^{-1}$, for $m=1,\dots, M$, and times $0=t_0<t_1<\dots< t_M$ such that if
    \begin{align}\label{cond:non:degen}
              \veps:=\inf_{1\leq m<M+1}\veps_m>0,
    \end{align}
and $\mu_m>0$, $N_m\geq1$ satisfy
    \begin{align}\label{cond:mu:N:converge}
        \gam_{m}\leq \frac{\mu_m}{\nu_0}\leq\gam N_m^2,
    \end{align}
for each $m=1,\dots, M$, then $\inf_{1\leq m<M+1}\nu_m>0$, and
    \begin{align}\label{eq:converge}
        |\nu_{m+1}-\nu|\leq\be|\nu_m-\nu|,
    \end{align}
for some $\be\in(0,1)$, for all $m=0,\dots,M-1$, where $\nu_m$ is given by \eqref{def:nu:sequence}.
\end{Thm}

Operationally, the ``non-degeneracy conditon'' \eqref{cond:non:degen} can be checked only after having chosen $\mu_m$. Once $\mu_m$ has been chosen, then the lower bound in the ``tuning condition'' \eqref{cond:mu:N:converge} can be checked. If $\mu_m$ does not satisfy the lower bound in \eqref{cond:mu:N:converge}, then there are three reasonable options: either wait longer to evaluate $\veps_m$ with the hope of obtaining a more suitable value for $\gam_m$, change the value of $\mu_m$ and repeat the process, or else collect observations on additional modes. In practice, it often suffices to simply wait for a judicious moment to evaluate $\veps_m$ and apply \eqref{def:new:rep}.

We will also prove another result for the time-averaged update formula. To state it, we define
    \begin{align}\label{def:bar:eps}
        \bar{\veps}_m:=\left|\frac{1}{t_{m}-t_{m}'}\int_{t_{m}'}^{t_{m}}\lb AP_n\tu(s),P_N(\tu(s)-u(s))\rb ds\right|\bar{\nu}_0^{-2},
    \end{align}
where $t_m>t_m'\geq t_{m-1}$, and suppose that each $\nu_m=\bar{\nu}_m$ is related as follows
    \begin{align}\label{def:nu:sequence:second}
      \bar{\nu}_{m+1}=\bar{\nu}_m+\mu\frac{\frac{1}{t_{m+1}-t_{m+1}'}\int_{t_{m+1}'}^{t_{m+1}}\Sob{P_Nw(s)}{L^2}^2ds}{\frac{1}{t_{m+1}-t_{m+1}'}\int_{t_{m+1}'}^{t_{m+1}}\lb P_N(-\De)\tu(s),P_Nw(s)\rb_{L^2} ds}.
    \end{align}

\begin{comment}
\begin{Thm}\label{thm:converge:second}
Given any $\nu>0$ and $M\geq1$, there exists a tuning parameter $\mu_0>0$, and a sequence of averaging windows, $J_m=[t_m',t_m)$, where $t_{m+1}'>t_m>0$, for $m=1,\dots, M-1$, such that whenever
            \begin{align}\label{cond:non:degen:second}
              \inf_{1\leq m<M+1}\bar{\veps}_m>0,
             \end{align}
for some $\veps>0$, then there exists $\be\in(0,1)$ such that for all $\mu\geq \mu_0$, there exists an observational threshold $N_0\geq1$ such that for any $N\geq N_0$, one has
    \begin{align}\label{eq:converge:second}
        |\bar{\nu}_{m+1}-\nu|\leq\be|\bar{\nu}_m-\nu|,
    \end{align}
for all $m=0,\dots,M-1$, where $\bar{\nu}_m$ is given by \eqref{def:mod:rep}. Moreover, $\{\bar{\nu}_m\}_{m=0}^M$ is a positive sequence of numbers.
\end{Thm}
\end{comment}

\begin{Thm}\label{thm:converge:second}
Assume the hypotheses of \cref{thm:ng:bounds}. Given any $\bar{\nu}>0$ and $M\geq1$, there exist constants $\bar{\gam},\bar{\gam}_1,\dots,\bar{\gam}_M>0$, where $\bar{\gam}$ depends only on the initial relative error $\de_0:=|\bar{\nu}_0-\nu|\nu^{-1}$ and  Grashof-type number, $|g|\bar{\nu}_0^{-2}$, and where $\bar{\gam}_m$ depends additionally on $\bar{\veps}_m^{-1}$ and non-dimesionalized frequency $(\bar{\nu}_0|J_m|)^{-1}$, for $m=1,\dots, M$, times $0\leq t_{m-1}\leq t_m'<t_m$, for all $m=1,\dots, M$, and averaging windows $[t_m',t_m)\subset[t_{m-1},t_m)$, for $m=1,\dots, M$, such that if
            \begin{align}\label{cond:non:degen:second}
              \inf_{1\leq m<M+1}\bar{\veps}_m>0,
             \end{align}
and $\bar{\mu}_m>0, N_m\geq1$ satisfy
    \begin{align}\label{cond:mu:N:converge:second}
        \bar{\gam}_m\leq\frac{\bar{\mu}_m}{\bar{\nu}_0}\leq\bar{\gam}N_m^2,
    \end{align}
for each $m=0,\dots,M-1$, then $\inf_{1\leq m<M+1}\bar{\nu}_m>0$, and 
    \begin{align}\label{eq:converge:second}
        |\bar{\nu}_{m+1}-\nu|\leq\be|\bar{\nu}_m-\nu|,
    \end{align}
for some $\be\in(0,1)$, for all $m=0,\dots, M-1$, where $\bar{\nu}_m$ is given by \eqref{def:nu:sequence:second}.
\end{Thm}

\begin{Rmk}\label{rmk:uniform}
Note that the conditions \eqref{cond:mu:N:converge} can be made uniform in $m$. Indeed, one may take $\gam':=\sup_{1\leq m<M+1}\gam_m$ and deduce the same conclusion \eqref{eq:converge} provided that $\mu,N>0$ satisfy
    \begin{align}\label{cond:mu:N:converge:uniform}
        \gam'\leq\frac{\mu}{\nu_0}\leq\gam N^2.
    \end{align}
In this case, $\gam'$ can be seen to depend on $\veps^{-1}$. A similar modification can be made for condition \eqref{cond:mu:N:converge:second}.
\end{Rmk} 

We obtain the following convergence results as immediate corollaries .%As an immediate corollary to \cref{thm:sensitivity} and \cref{thm:converge}, we have the following convergence result.

\begin{Cor}\label{cor:converge}
Under the hypothesis of \cref{thm:converge}, if $M=\infty$ and \eqref{cond:non:degen} holds, then $\{\nu_m\}_{m=0}^\infty$ converges exponentially to $\nu$. Similarly, under the hypothesis of \cref{thm:converge:second}, if $M=\infty$ and \eqref{cond:non:degen:second} holds, then $\{\bar{\nu}_m\}_{m=0}^\infty$ converges exponentially to $\nu$.
\end{Cor}

The proofs of \cref{thm:converge} and \cref{thm:converge:second} are carried out in \cref{sect:proof}.

\begin{Rmk}\label{rmk:algorithm}
Five remarks regarding the practicality of the hypotheses of \cref{thm:converge} are in order. Firstly, in numerical experiments \cite{CarlsonHudsonLarios2020, CarlsonHudsonLariosMartinezNgWhitehead2021} when the initial guess for viscosity was taken larger than its true value it was observed that such a choice typically leads to better-behaved convergence. The proof of \cref{thm:converge} reflects this choice. Interestingly, the assumption is not a matter of convenience in the proof and is important to make in order to guarantee, in general, that the first update produces a valid (positive) and improved approximation to the true value.

Secondly, the restrictions on $\mu$ depend on the true value of the viscosity, which is not known apriori. However, it is often presumed that it is known to belong to some neighborhood of values. One purpose for the parameter $\al$ that appears in the hypothesis \cref{thm:sensitivity} is to allow for enough flexibility in estimates to account for this. Also, in numerical experiments, $\mu$ is often taken to be quite large relative to the reference viscosity. Although the theorem suggests that larger values of $\mu$ require a larger number of observations, the relationship suggested by the theorem is often pessimistic. Indeed, \cite{GeshoOlsonTiti2016, FarhatGlattHoltzMartinezMcQuarrieWhitehead2020} suggest that far less observations are needed in practice that what is rigorously required (see also the fifth remark below).

Thirdly, although we cannot theoretically guarantee that the non-degeneracy condition \eqref{cond:non:degen} holds, it can nevertheless be practically implemented. In particular, at the proposed time of update, the condition can be checked, at which time a suitable adjustment for $\mu$ can be made so that \eqref{eq:converge} still holds (see \eqref{cond:mu:proof}) or else, one may wait longer. In practice, one can set a rule so that the update is only applied under favorable conditions (see \cite{CarlsonHudsonLariosMartinezNgWhitehead2021}) From this point of view, the condition \eqref{cond:non:degen} is a practical one.

Fourthly, the size of $\mu$ depends inversely on the non-degeneracy parameter, $\veps$. On the other hand, the relation \eqref{cond:mu:N:nu} is always maintained. Thus, if $\veps\goesto0$, then $\mu\goesto\infty$ and subsequently $N\goesto\infty$. From the point of view of ill-posedness, the problem of asymptotically inferring $\nu$ from the observations $\{P_Nu(t)\}_{t\geq0}$ becomes ``more ill-posed'' as $\veps\goesto0$. Thus, one interpretation of the chain of dependencies between $\mu, N,\veps$ is that the increasing degeneracy of the problem can be balanced by having access to observations of higher-resolution.

Lastly, the size of the constants $\gam_m,\bar{\gam}_m$ and $\gam,\bar{\gam}$, which dictate how $\mu$ should be tuned and ultimately, how many observations one needs, can be tracked down in the analysis. Essentially, $\bar{\gam}_m\gtrsim G^2$. In the context of turbulence, it has been shown that $G\sim \text{Re}^2$, where $Re$ is the Reynolds number of the flow (see \cite{DascaliucFoiasJolly2008, DascaliucFoiasJolly2009}), so one may expect the lower bounds for $\bar{\gam}_m$ to be quite large in general. Indeed, in the context of two-dimensional turbulence, the dissipative cut-off of the enstrophy spectrum is posited to occur at wave-numbers of size on the order of $G^{1/4}$ \cite{DascaliucFoiasJolly2008}. Thus, the number of the observations required \eqref{cond:mu:N:converge} or \eqref{cond:mu:N:converge:second} may seem impratically large. However, the rigourous analysis is always representative of the ``worst-case scenario."  On the other hand, it is reasonable to expect that such worst-case scenarios are not typically realized in practice; support for this hypothesis initially was demonstrated in the work \cite{GeshoOlsonTiti2016}. The numerical studies carried out in \cite{CarlsonHudsonLarios2020} for the nudging-based corroborated the observations in \cite{GeshoOlsonTiti2016}, but in the context of the nudging-based approach to parameter estimation studied here by observing that convergence nevertheless occurs in regimes that were much less pessimistic than the ones suggested by bounds identified here.
\end{Rmk}

\begin{Rmk}\label{rmk:real}
Although the setting considered here is quite ideal, certain elements can be made more physical. One such element is the setting of non-spectral observations, such as local spatial averages, or nodal value measurements of the velocity field. We believe that the approach taken here can also be adapted to accommodate such observations, though it will require a significant technical effort.

Another physical aspect that one may try to incorporate is the presence of noise in the observations. This was studied from the point of view of data assimilation in \cite{BlomkerLawStuartZygalakis2013, BessaihOlsonTiti2015} and from the point of view of parameter estimation in \cite{CialencoGlattHoltz2011}. However, due to the need for estimating time-derivatives of \eqref{eq:nse:ff}, \eqref{eq:nse:ng:ff}, the analysis required to establish \cref{thm:converge} does not appear to be adaptable in a straight-forward manner. From this point of view, the update rule given by \eqref{def:mod:rep} appears to be more amenable to accommodating the situation of noisy observations. Further investigation on these matters are reserved for a future work.
\end{Rmk}

\begin{Rmk}\label{rmk:convention}
For the remainder of the article, we will suppress the notation $\T^2$ when writing Sobolev spaces. It is understood that we are working in the domain $[0,2\pi]^2$, where \eqref{eq:nse:ff} and \eqref{eq:nse:ng:ff} are supplemented with periodic boundary conditions.
\end{Rmk}

\section{Convergence Analysis}\label{sect:proof}

The main goal of this section is to prove \cref{thm:converge} and \cref{thm:converge:second}. In preparation for this, we will reformulate the theorem with notation that will be convenient for its proof. In the process, we also motivate the main idea behind the formulas \eqref{def:new:rep} and \eqref{def:mod:rep}.

\subsection{Derivation of the update formula}\label{sect:derive}
Let $\nu,\tnu>0$ and $g\in H_\s^1$.  Given $\mu>0$ and $N\geq1$ satisfying \eqref{cond:mu:N:nu} and \eqref{cond:mu:N:tnu}, let $u,\tu$ denote smooth solutions to \eqref{eq:nse:ff}, \eqref{eq:nse:ng:ff}, respectively. Let $w:=\tu-u$ denote the ``sensitivity variable". Then
    \begin{align}\notag
        \frac{dw}{dt}+\tnu A\tu-\nu Au+B(w,w)+(DBu)w&=-\mu P_Nw,\quad w(0)=w_0:=\tu_0-u_0,
    \end{align}
where
    \begin{align}\label{def:DB}
        (DBv_1)v_2=B(v_1,v_2)+B(v_2,v_1).
    \end{align}
Equivalently, the sensitivity variable satisfies
    \begin{align}\label{eq:nse:diff:ng}
        \frac{dw}{dt}+\nu Aw+B(w,w)+(DBu)w&=-(\De\nu) A\tu-\mu P_Nw,
    \end{align}
where
    \begin{align}\label{def:diff:nu}
        \De\nu=\tnu-\nu.
    \end{align}
Applying the low-pass filter, $P_N$, to \eqref{eq:nse:diff:ng}, we obtain
    \begin{align}\label{eq:nse:diff:low}
        \frac{dw_N}{dt}+\nu Aw_N+B_N(w,w)+(DB_Nu)w=-(\De\nu)A\tu_N-\mu w_N,\quad w_N(0)=P_Nw_0,
    \end{align}
where $w_N=P_Nw$ and $B_N=P_NB$. Upon taking the scalar product of \eqref{eq:nse:diff:low} with $w_N$ and then $Aw_N$, one obtains the scalar balance equations
    \begin{align}\label{eq:diff:balance:N}
        \begin{split}
        \frac{1}2\frac{d}{dt}|w_N|^2+\nu\lVert w_N\rVert^2+\mu|w_N|^2&=-\lb B_N(w,w),w_N\rb-\lb (DB_Nu)w,w_N\rb\\
        &\quad-(\De\nu)\lb A\tu_N, w_N\rb,\\
        \frac{1}2\frac{d}{dt}\lVert w_N\rVert^2 +\nu|Aw_N|^2+\mu\lVert w_N\rVert^2&=-\lb B_N(w,w),Aw_N\rb-\lb (DB_Nu)w,Aw_N\rb\\
        &\quad-(\De\nu)\lb A\tu_N,Aw_N\rb.
        \end{split}
    \end{align}
We introduce the following notation
    \begin{align}\label{def:diff:energy}
        \begin{split}
        \EE(t):=\frac{1}2|w(t)|^2,\hspace{10pt}\quad\EE_N(t)&:=\frac{1}2|w_N(t)|^2,\quad\ \ \ \EEd_N(t):=\frac{d}{dt}\EE_N(t)\\
        \ZZ(t):=\frac{1}2\lVert w(t)\rVert^2,\hspace{4pt}\quad\ZZ_N(t)&:=\frac{1}2\lVert w_N(t)\rVert^2,\quad\ \ZZd_N(t):=\frac{d}{dt}\ZZ_N(t)\\
        \PP(t):=\frac{1}2| Aw(t)|^2,\quad\PP_N(t)&:=\frac{1}2| Aw_N(t)|^2,\quad \PPd_N(t):=\frac{d}{dt}\PP_N(t),
        \end{split}
    \end{align}
and
 \begin{align}\label{def:J}
        \begin{split}
        \JJ_1&:=\lb B(w,w_N),Q_Nw\rb-\lb (DB_Nu)w,w_N\rb-(\De\nu)\lb A\tu_N, w_N\rb\\
        \JJ_2&:=-\lb B_N(w,w),Aw_N\rb-\lb (DB_Nu)w,Aw_N\rb-(\De\nu)\lb A\tu_N, Aw_N\rb.
        \end{split}
    \end{align}
Then, upon integrating by parts, we see that \eqref{eq:diff:balance:N} may be rewritten as
    \begin{align}\label{eq:diff:balance:Ndot}
        \begin{split}
       \EEd_N+2\mu\EE_N+2\nu\ZZ_N&=\JJ_1,\\
       \ZZd_N+2\mu\ZZ_N+2\nu\PP_N&=\JJ_2.
       \end{split}
    \end{align}
   
Retaining only the terms that depend at most linearly on the difference, $w$, as well as the feedback control term, $2\mu\EE_N$, one obtains
    \begin{align}\notag
        2\mu\EE_N\approx-(\De\nu)\lb A\tu, w_N\rb.
    \end{align}
Upon solving for $\nu$, one thus derives an approximation for the true viscosity in terms of $\tnu$ given by
    \begin{align}\label{eq:approx:first}
        \nu\approx \tnu+2\mu\frac{\EE_N}{\lb A\tu_N, w_N\rb}.
    \end{align}
Based on this, an algorithm for successively providing approximations to the true viscosity can be developed.

Given a sequence of viscosities, $\nu_0,\nu_1,\dots, \nu_M>0$, nudging parameters, $\mu_1,\dots, \mu_M>0$, observations densities, $N_1,\dots, N_M\geq1$, and times $0=t_0<\dots< t_M$, where $M\geq1$, let $t_m$ denote the $m$--th ``update time" and $I_m=[t_m,t_{m+1})$, where $m=0,\dots, M-1$. In the framework of \eqref{eq:nse:ng:I0}, \eqref{eq:nse:ng:Im}, we initially consider, for $t\in I_0$, the system 
    \begin{align}\label{eq:diff:balance:Ndot:I0}
        \begin{split}
       \EEd_N+2\mu_1\EE_N+2\nu\ZZ_N&=\JJ_1^{(0)},\quad \EE_N(0)=\frac{1}2|w_0|^2,\\
       \ZZd_N+2\mu_1\ZZ_N+2\nu\PP_N&=\JJ_2^{(0)},\quad \ZZ_N(0)=\frac{1}2\lVert w_0\rVert^2,
       \end{split}
    \end{align}
where the dependence on $\De\nu$ in $\JJ_1^{(0)},\JJ_2^{(0)}$ is replaced by $\De\nu_0=\nu_0-\nu$. Then for $t\in I_m$, where $m\geq1$, we consider
    \begin{align}\label{eq:diff:balance:Ndot:Im}
        \begin{split}
       \EEd_N+2\mu_{m+1}\EE_N+2\nu\ZZ_N&=\JJ_1^{(m)},\quad \EE_N(t_m)=\frac{1}2|w_N^{(m)}|^2,\\
       \ZZd_N+2\mu_{m+1}\ZZ_N+2\nu\PP_N&=\JJ_2^{(m)},\quad \ZZ_N(t_m)=\frac{1}2\lVert w_N^{(m)}\rVert^2,
       \end{split}
    \end{align}
where $\De\nu$ in $\JJ_1^{(m)},\JJ_2^{(m)}$ is given by $\De\nu^m=\nu^m-\nu$, and we let
    \begin{align}\label{def:Nstep:m}
        w_N^{(m)}=w_N(t_{m}^-),\quad \tu_N^{(m)}=\tu_N(t_{m}^-).
    \end{align}

\subsubsection{Update rule $\#1$}
With this notation in hand, the update formula \eqref{def:new:rep} is therefore given in a recursive form by
    \begin{align}\label{def:update}
        \nu_{m+1}=\nu_m+2\mu\frac{\EE_N^{(m+1)}}{\lb A\tu_N^{(m+1)},w_N^{(m+1)}\rb},
    \end{align} 
for $m=0,\dots, M-1$. Ultimately, for the above considerations to be well-defined and provide increasingly accurate approximations to $\nu$, we  will establish two properties of \eqref{def:update}:
    \begin{enumerate}
        \item the sequence $\{\nu_m\}_{m=0}^M$ is positive, so that \eqref{eq:diff:balance:Ndot:Im} is well-defined over each time interval $I_m$;
        \item the error between successive updates decrement in a geometric fashion.
    \end{enumerate}
We will show that the first property holds in the process of establishing the second property. The proof will ultimately proceed by induction.

\subsubsection{Update rule $\#2$}
Upon taking the time average, instead, of the first equation in \eqref{eq:diff:balance:Ndot} over the interval $[t',t]$, where $t-t'\gg1$, then retaining only terms that are at most linear in $w$, as well as the term with pre-factor $\mu$, one may then solve for $\nu$ to obtain
    \begin{align}\label{eq:approx:second}
        \nu\approx\tnu+2\mu\frac{\frac{1}{t-t'}\int_{t'}^t\EE_N(s)ds}{\frac{1}{t-t'}\int_{t'}^t\lb A\tu_N(s),w_N(s)\rb ds}.
    \end{align}
Hence, in the context of \eqref{eq:diff:balance:Ndot:I0}--\eqref{def:Nstep:m}, the second update rule can be written as
    \begin{align}\label{def:update:second}
        \bar{\nu}_{m+1}=\bar{\nu}_m+2\mu\frac{\frac{1}{|J_{m+1}|}\int_{J_{m+1}}\EE_N(s)ds}{\frac{1}{|J_{m+1}|}\int_{J_{m+1}}\lb A\tu_N(s),w_N(s)\rb ds},
    \end{align}
where $J_m=[t_m',t_m)$ such that $t_{m+1}'\geq t_m>t_m'$, for all $m\geq1$. Recall that $I_m=[t_m,t_{m+1})$, for all $m\geq0$, where $t_0=0$. Thus $J_m\subset I_{m-1}$, for all $m\geq1$. As in the first update rule, we must show that each newly proposed viscosity is admissiable, i.e., positive. The proof will also proceed by induction, though it will be far simpler due to fact the averaging process dispenses with the need for careful time-derivative estimates.

\subsection{Outline of the convergence argument}
In this section, we outline the convergence arguments for each update rule to provide a sense for the expected difficulties therein.

\subsubsection{Update rule $\#1$}
Observe that we may solve for $\De\nu$ in \eqref{eq:diff:balance:Ndot} to obtain
    \begin{align}\label{eq:diff:nu}
        \De\nu=\frac{-\EEd_N-2\nu\ZZ_N+\lb B(w,w_N),Q_Nw\rb-\lb (DB_Nu)w,w_N\rb-2\mu\EE_N}{\lb A\tu_N,w_N\rb}.
    \end{align}
Considering \eqref{eq:diff:nu} over the interval $I_m$, then evaluating at $t_{m+1}^-$, we obtain
    \begin{align}\label{eq:diff:nu:n}
        &\nu_{m}-\nu\\
        &={\scalebox{1.17}{$\frac{-\EEd_N^{(m+1)}-2\nu\ZZ_N^{(m+1)}+\lb B(w^{(m+1)},w_N^{(m+1)}),Q_Nw^{(m+1)})\rb-\lb (DB_Nu^{(m+1)})w^{(m+1)},w_N^{(m+1)}\rb-2\mu\EE_N^{m+1}}{\lb A\tu_N^{(m+1)},w_N^{(m+1)}\rb}$}},\notag
    \end{align}
where
    \begin{align}\label{def:step:m}
        u^{(m)}=u(t_{m}^-),\quad \EEd_N^{(m)}=\EEd_N(t_{m}^-).
    \end{align}
Assuming for the moment that the viscosities defined by \eqref{def:update} for $m=1,\dots, M-1$ are indeed positive, we may then combine \eqref{def:update} and \eqref{eq:diff:nu:n} to obtain
    \begin{align}
        &\nu_{m+1}-\nu=(\nu_{m+1}-\nu_m)+(\nu_m-\nu)\notag\\
        &=\frac{-\EEd_N^{(m+1)}-2\nu\ZZ_N^{(m+1)}+\lb B(w^{{(m+1)}},w_N^{{(m+1)}}),Q_Nw^{{(m+1)}}\rb-\lb (DB_Nu^{{(m+1)}})w^{{(m+1)}},w_N^{(m+1)}\rb}{\lb A\tu_N^{(m+1)},w_N^{(m+1)}\rb}.\notag
    \end{align}
For $M\geq1$, we suppose that $\veps_m:=\nu_0^{-2}|\lb A\tu_N^{(m)},w_N^{(m)}\rb|>0$. Then
    \begin{align}\label{eq:induction:step}
        &|\nu_{m+1}-\nu|\leq \frac{1}{\veps_{m+1}\nu_0^2}\left(\left|\EEd_N^{(m+1)}\right|+2\nu\ZZ_N^{(m+1)}+\left|\lb B(w^{(m+1)},w_N^{(m+1)}),Q_Nw^{(m+1)}\rb\right|\right.\notag\\
        &\qquad\qquad\qquad\qquad\left.+\left|\lb (DB_Nu^{(m+1)})w^{(m+1)},w_N^{(m+1)}\rb\right|\right).
    \end{align}
We will ultimately show that for each $m=1,\dots, M-1$, the right-hand side can be controlled so that
    \begin{align}\label{eq:desire}
        |\nu_{m+1}-\nu|\leq C\veps_{m+1}^{-1}\nu_0\left(\frac{\nu_0}{\mu_m}\right)^p|\nu_m-\nu|,
    \end{align}
for some $p\in(0,1)$ and some constant $C>0$, independent of $m$, provided that $t_{m+1}$, $\mu>0$, and $N\geq1$ are chosen sufficiently large. We will then deduce the desired error estimate by ultimately taking $\mu$ large enough so that $\mu\geq \nu_0 C^{1/p}(\be\veps_{m+1})^{-1/p}$, for some $\be$ sufficiently small.

The term $\EEd_N^{(m+1)}$ constitutes the heart of the difficulty of the proof. A significant effort of the proof \cref{thm:converge} is dedicated to this task alone and is the content of \cref{lem:power} below.

\subsubsection{Update rule $\#2$}
We recall that $J_m\subset I_{m-1}$, for all $m\geq1$. Thus, we obtain the corresponding analog of \eqref{eq:diff:nu:n} by considering \eqref{eq:diff:balance:Ndot} over $I_m$, then taking the time-average over $J_{m+1}$, then solving for $\De$ to obtain
    \begin{align}
     \bar{\nu}_{m}-\nu
        &={\scalebox{1.17}{$-\frac{\frac{\EE_N(t_{m+1})-\EE_N(t_{m+1}')}{|J_{m+1}|}}{\frac{1}{|J_{m+1}|}\int_{J_{m+1}}\lb A\tu_N(s),w_N(s)\rb ds}$}}\notag\\
        &\ \quad{\scalebox{1.17}{$+\frac{\frac{1}{|J_{m+1}|}\int_{J_{m+1}}-2\nu\ZZ_N(s)+\lb B(w(s),w_N(s)),Q_Nw(s))\rb-\lb (DB_Nu(s))w(s),w_N(s)\rb-2\mu\EE_N(s) ds}{\frac{1}{|J_{m+1}|}\int_{J_{m+1}}\lb A\tu_N(s),w_N(s)\rb ds}$}}.\notag
    \end{align}
Then
    \begin{align}
        \bar{\nu}_{m+1}-\nu&=(\nu_{m+1}-\nu_m)+(\nu_m-\nu)\notag\\
        &=-\frac{\frac{\EE_N(t_{m+1})-\EE_N(t_{m+1}')}{|J_{m+1}|}}{\frac{1}{|J_{m+1}|}\int_{J_{m+1}}\lb A\tu_N(s),w_N(s)\rb ds}\notag\\
        &\quad{\scalebox{1.17}{$+\frac{\frac{1}{|J_{m+1}|}\int_{J_{m+1}}-2\nu\ZZ_N(s)+\lb B(w(s),w_N(s)),Q_Nw(s))\rb-\lb (DB_Nu(s))w(s),w_N(s)\rb ds}{\frac{1}{|J_{m+1}|}\int_{J_{m+1}}\lb A\tu_N(s),w_N(s)\rb ds}$}}.\notag
    \end{align}
Hence, if we suppose that $\bar{\veps}_{m+1}:=\bar{\nu}_0^{-2}\left|\frac{1}{|J_{m+1}|}\int_{J_{m+1}}\lb A\tu_N(s),w_N(s)\rb ds\right|$, then
    \begin{align}\label{eq:induction:step:second}
        |\bar{\nu}_{m+1}-\nu|\leq& \frac{1}{\bar{\veps}_{m+1}\bar{\nu}_0^2}\left(\frac{|\EE_N(t_{m+1})-\EE_N(t_{m+1}')|}{|J_{m+1}|}\right.\\
        &\left.{\scalebox{0.9}{$+\frac{1}{|J_{m+1}|}\int_{J_{m+1}}2\nu\ZZ_N(s)+|\lb B(w(s),w_N(s)),Q_Nw(s))\rb|+|\lb (DB_Nu(s))w(s),w_N(s)\rb| ds$}}\right).\notag
    \end{align}
Notice that in comparison with \eqref{eq:induction:step}, no time-derivative term appears in \eqref{eq:induction:step:second}. This particular feature will greatly simplify the proof of convergence. In exchange for not having to estimate time-derivatives, we will instead impose a uniform lower bound on the size of the averaging windows, $|J_{m+1}|$.

\subsection{Estimates for the power of the error on low frequencies} 
This section is dedicated to establishing estimates on the \textit{power} generated by the sensitivity variable, $w$, at small frequencies, $|k|\leq N$, that is, we obtain estimates on the quantity $\EEd_N$. This will be done by identifying a dissipative structure for its square.

To state this lemma precisely, we must make a subtle distinction in the defining of the evolution equation that governs $\EEd_N^2$: Given $\nu_0,\nu_1,\dots,\nu_M>0$, and $0=t_0<\dots<t_M$, where $M\geq1$, we consider $\EE_N,\ZZ_N,\PP_N$ such that \eqref{eq:diff:balance:Ndot:I0} holds when $t\in I_0$ and \eqref{eq:diff:balance:Ndot:Im} holds when $t\in I_m$, for $m\geq1$. It follows from \eqref{eq:diff:balance:Ndot} that the evolution of $\EEd_N^2$ over $t\in I_m$ is determined by
     \begin{align}\notag
       \frac{d}{dt}\EEd_N^2+2\mu\frac{d}{dt}(\EE_N\EEd_N)+2\nu\frac{d}{dt}(\ZZ_N\EEd_N)=\JJd_1\EEd_N+\JJ_1\EEdd_N.
    \end{align}
Further expanding with \eqref{eq:diff:balance:Ndot} yields the evolution equation
    \begin{align}\label{eq:diff:balance:Ndot:squared}
        \frac{d}{dt}\EEd_N^2+2\mu\EEd_N^2+2\nu\ZZd_N\EEd_N+2\left(\mu\EE_N+\nu\ZZ_N\right)\EEdd_N=\JJd_1\EEd_N+\JJ_1\EEdd_N.
    \end{align}
Note that the ``initial condition" is determined by evaluating \eqref{eq:diff:balance:Ndot} at $t=t_m^+$, which is well-defined since we will always have $u(t_m),\tu(t_m)\in H^3_\s$, for each $m$.

%\EEd_N^2(t_m):=\left(\frac{d}{dt}\bigg{|}_{t_m^+}\EE_N\right)^2.

\begin{Lem}\label{lem:power}
Let $\tau_0>0$ and $g\in H_\s^2(\T^2)$. Suppose that $u(\tau_0)\in B_\ell(R_\ell)$ and $\tu(\tau_0)\in B_\ell(\al R_\ell)$, for $\ell=1,2,3$, for some $\al\geq1$, where $R_k$ is given by \eqref{def:rad:Hk}. Suppose that $\mu, N,\nu,\tnu$ satisfy 
    \begin{align}
    \eqref{cond:mu:N:tnu}, \eqref{cond:mu:N:nu}, \eqref{cond:mu:ng:H2}, \eqref{cond:mu:ng:Hk}_3, \eqref{cond:mu:sensitivity1}, \eqref{cond:mu:sensitivity2}.\label{cond:mu:N:nu:tnu}
    \end{align}
Let $u(\cdotp;u(\tau_0), \nu)$ and $\tu(\cdotp;\tu(\tau_0),\tnu)$ denote the corresponding global-in-time unique solutions of \eqref{eq:nse:ff} and \eqref{eq:nse:ng:ff} over $[\tau_0,\infty)$ corresponding to initial data $u(\tau_0),\tu(\tau_0)$ and viscosities $\nu,\tnu$, respectively. Suppose furthermore that $\nu,\tnu$ satisfy
    \begin{align}\label{cond:nu:error}
        \frac{|\De\nu|}{\nu}\leq \de,
    \end{align}
for some $\de>0$, where $\De\nu$ is given by \eqref{def:diff:nu}. There exists a constant $c_4\geq1$ and time $\tau_1\geq\tau_0$ such that if
    \begin{align}\label{cond:mu:power}
        \mu\geq c_4(G\vee1)^2\nu,
    \end{align}
and either 
    \begin{align}\label{cond:diff:nu}
        \de\leq 1\quad\text{or}\quad \De\nu\geq0,
    \end{align}
then
    \begin{align}\notag
        |\EEd_N(t)|\leq \nu^3\left(\frac{\nu}{\mu}\right)^{1/2}\left(\frac{|\De\nu|}{\nu}\right)L,
    \end{align}
for all $t\geq\tau_1$, where $L$ is some constant that depends on $\al, \nu\tnu^{-1}, \de, G, \tG$.
\end{Lem}

\begin{Rmk}\label{rmk:lorenz}
A result similar to \cref{lem:power} was developed in \cite{CarlsonHudsonLariosMartinezNgWhitehead2021}, where parameter estimation for Lorenz 63 system was studied. Indeed, in this case, convergence of the algorithm reduced to establishing sufficient control on the parameter error thorugh a relation similar to \eqref{eq:induction:step}. The key step in \cite{CarlsonHudsonLariosMartinezNgWhitehead2021} was to obtain estimates for the time-derivative of the state through a Lyapunov-type functional, the evolution of which could be developed in a way similar to \eqref{eq:diff:balance:Ndot:squared}; the desired estimates were deduced suitably from there, as we are about to do now for the case of the 2D NSE. In addition to being technically more involved, one may identify one subtle difference between their proofs, namely, the size of the initial guess for the viscosity plays a distinguished role here. We refer the reader back to \cref{rmk:algorithm} to recall the discussion regarding this point.
\end{Rmk}

Before proceeding to the proof of \cref{lem:power}, we will first expand \eqref{eq:diff:balance:Ndot:squared} to separate the sign-definite terms from the sign-indefinite terms. To this end, let us first develop the left-hand side.

First observe that
    \begin{align}\label{eq:Ed:Zd}
        \begin{split}
        2\nu\ZZd_N\EEd_N&=2\nu\left(-2\mu\ZZ_N-2\nu\PP_N+\JJ_2\right)\left(-2\mu\EE_N-2\nu\ZZ_N+\JJ_1\right)\\
        &=8\mu^2\nu\EE_N\ZZ_N+8\mu\nu^2\left(\ZZ_N^2+\EE_N\PP_N\right)+8\nu^3\ZZ_N\PP_N\\
        &\quad-4\nu\left(\mu\ZZ_N+\nu\PP_N\right)\JJ_1-4\nu(\mu\EE_N+\nu\ZZ_N)\JJ_2+2\nu\JJ_1\JJ_2,
        \end{split}
    \end{align}
where $\JJ_1$ and $\JJ_2$ are defined in \eqref{def:J}. Next observe that by differentiating \eqref{eq:diff:balance:Ndot}, we obtain
    \begin{align}\label{eq:diff:balance:Nddot}
       \EEdd_N&=-2\mu\EEd_N-2\nu\ZZd_N+\JJd_1\notag\\
       &=\left(4\mu^2\EE_N+4\mu\nu\ZZ_N-2\mu\JJ_1\right)+\left(4\mu\nu\ZZ_N+4\nu^2\PP_N-2\nu\JJ_2\right)+\JJd_1\notag\\
       &=4\mu^2\EE_N+8\mu\nu\ZZ_N+4\nu^2\PP_N-2\mu\JJ_1-2\nu\JJ_2+\JJd_1.
    \end{align}
It follows that
    \begin{align}
        \EE_N\EEdd_N&=4\mu^2\EE_N^2+8\mu\nu\EE_N\ZZ_N+4\nu^2\EE_N\PP_N-2\mu\EE_N\JJ_1-2\nu\EE_N\JJ_2+\EE_N\JJd_1\notag\\
        \ZZ_N\EEdd_N&=4\mu^2\EE_N\ZZ_N+8\mu\nu\ZZ_N^2+4\nu^2\ZZ_N\PP_N-2\mu\ZZ_N\JJ_1-2\nu\ZZ_N\JJ_2+\ZZ_N\JJd_1.\notag
    \end{align}
In particular
    \begin{align}\label{eq:Ndot:diss}
        &\left(2\mu\EE_N+2\nu\ZZ_N\right)\EEdd_N\notag\\
        &=8\mu^3\EE_N^2+16\mu^2\nu\EE_N\ZZ_N+8\mu\nu^2\EE_N\PP_N-4\mu^2\EE_N\JJ_1-4\mu\nu\EE_N\JJ_2+2\mu\EE_N\JJd_1\notag\\
        &\quad+8\mu^2\nu\EE_N\ZZ_N+16\mu\nu^2\ZZ_N^2+8\nu^3\ZZ_N\PP_N-4\mu\nu\ZZ_N\JJ_1-4\nu^2\ZZ_N\JJ_2+2\nu\ZZ_N\JJd_1.
    \end{align}
Lastly, observe that
    \begin{align}\label{eq:JJd1:EEd}
    \JJd_1\EEd_N=-2\mu\EE_N\JJd_1-2\nu\ZZ_N\JJd_1+\JJ_1\JJd_1
    \end{align}
Finally, upon returning to \eqref{eq:diff:balance:Ndot:squared} and applying \eqref{eq:Ed:Zd}--\eqref{eq:JJd1:EEd}, we obtain
    \begin{align}
        &\frac{d}{dt}\EEd_N^2+2\mu\EEd_N^2+\left(8\mu^2\nu\EE_N\ZZ_N+8\mu\nu^2\ZZ_N^2+8\mu\nu^2\EE_N\PP_N+8\nu^3\ZZ_N\PP_N\right)\notag\\
        &\quad+\left[\left(8\mu^3\EE_N^2+16\mu^2\nu\EE_N\ZZ_N+8\mu\nu^2\EE_N\PP_N\right)+\left(8\mu^2\nu\EE_N\ZZ_N+16\mu\nu^2\ZZ_N^2+8\nu^3\ZZ_N\PP_N\right)\right]\notag\\
        &=4\nu\left(\mu\ZZ_N+\nu\PP_N\right)\JJ_1+4\nu(\mu\EE_N+\nu\ZZ_N)\JJ_2-2\nu\JJ_1\JJ_2\notag\\
        &\quad +\left(4\mu^2\EE_N\JJ_1+4\mu\nu\EE_N\JJ_2-2\mu\EE_N\JJd_1\right)+\left(4\mu\nu\ZZ_N\JJ_1+4\nu^2\ZZ_N\JJ_2-2\nu\ZZ_N\JJd_1\right)\notag\\
        &\quad -2\mu\EE_N\JJd_1-2\nu\ZZ_N\JJd_1+\JJ_1\JJd_1\notag\\
        &\quad+\left(4\mu^2\EE_N\JJ_1+8\mu\nu\ZZ_N\JJ_1+4\nu^2\PP_N\JJ_1-2\mu\JJ_1^2-2\nu\JJ_1\JJ_2+\JJ_1\JJd_1\right).\notag
    \end{align}
We let 
    \begin{align}\label{def:D}
        \DD:=8\mu^3\EE_N^2+32\mu^2\nu\EE_N\ZZ_N+16\mu\nu^2\EE_N\PP_N+16\nu^3\ZZ_N\PP_N+24\mu\nu^2\ZZ_N^2+2\mu\JJ_1^2,
    \end{align}
and
    \begin{align}\label{def:F}
        \FF&:=8\mu^2\EE_N\JJ_1+16\mu\nu\ZZ_N\JJ_1+8\nu^2\PP_N\JJ_1+8\mu\nu\EE_N\JJ_2+8\nu^2\ZZ_N\JJ_2\notag\\
        &\qquad-4\nu\JJ_1\JJ_2-4\mu\EE_N\JJd_1-4\nu\ZZ_N\JJd_1+2\JJ_1\JJd_1.
    \end{align}
From \eqref{def:D} and \eqref{def:F}, we may ultimately write
    \begin{align}\label{eq:diff:balance:Ndot:compact}
        \frac{d}{dt}\EEd_N^2+2\mu\EEd_N^2+\DD=\FF.
    \end{align}

The main goal, then, is to strike a tolerable balance between $\DD$ and $\FF$ that depends favorably in $\mu$ and $|\De\nu|$. In particular, it will suffice to show that $\FF$ can be bounded in terms of $|\De\nu|$, but \textit{independently} of $\mu$. The Gronwall inequality, will then yield bounds on $\EEd_N^2$ that depend proportionally on $|\De\nu|$, but inversely in $\mu$, as desired in \eqref{eq:desire}.  Specifically, let us consider
    \begin{align}\label{def:lam}
        \lam_1=\frac{\PP}{\nu^2},\quad \lam_2=\frac{|\De\nu|}{\nu},\quad \lam_3=\frac{\nu}{\mu}. 
    \end{align}
We claim that there exists a (generalized) tri-variate polynomial, $\QQ=\QQ(\lam_1,\lam_2,\lam_3)$, that is, whose exponents are positive real numbers, satisfying the following properties:
    \begin{itemize}
        \item $\QQ$ has no constant term;
        \item its coefficients depend only on the parameters $\s_\ell$, $G$, $\tnu\nu^{-1}$, $|\tu|\nu^{-1}$, $\lVert\tu\rVert\nu^{-1}$, $|A\tu|\nu^{-1}$, in an increasing fashion;
        \item $\lam_2$ always occurs jointly with either $\lam_1$ or $\lam_3$, i.e., $\lam_2$ never appears in isolation, as a monomial in itself,
        \item the degree of $\lam_1$ in any monomial is always $\geq1/2$,
        \item the degree of $\lam_2$ in any monomial is always $\geq1$,
    \end{itemize}
and there exist polynomials $\RR=\RR(\lam_2)$ and $\SS=\SS(\lam_3)$ satisfying:
    \begin{itemize}
        \item $\RR$, $\SS$ have no constant term;
        \item the coefficients of $\RR$ are constants that are independent of all material parameters;
        \item the coefficients of $\SS$ depend increasingly on $G$,
    \end{itemize}
such that
    \begin{align}\label{eq:Ndot:structure}
        \frac{d}{dt}\EEd_N^2+2\mu\EEd_N^2\leq \nu^7\QQ(\lam_1,\lam_2,\lam_3)+\left(\RR(\lam_2)+\SS(\lam_3)-1\right)\DD.
    \end{align}

Since we must deal with a plethora of terms, as captured by $\FF$ in \eqref{def:F}, in what follows, we will adopt the convention that $\QQ, \RR, \SS$ denote generic polynomials of the forms, respectively, as described above. They will behave like ``constants" in that they may change line-to-line in the estimates below after undergoing applications of the Cauchy-Schwarz inequality or Young's inequality, for instance. We will also make use of the following notation
    \begin{align}\label{def:tQQ}
        \tQQ=\QQ_0\frac{\PP}{\nu^2}\quad \text{or}\quad \QQ_0\frac{\nu}{\mu}\left(\frac{|\De\nu|}{\nu}\right)^2,\quad\text{for some}\ \QQ_0,
    \end{align}
where $\QQ_0$ is a polynomial of the form described above, except possibly independent of $\lam_1$, in the first case, or $\lam_2$ or $\lam_3$, in the second case.

Before proceeding to the proof, let us collect a few useful inequalities in terms of the notation developed above that demonstrate our convention. These all follow from interpolation:
    \begin{align}\label{est:diff:interpolation}
        \begin{split}
        \Sob{w}{L^4}&\leq C\ZZ^{1/4}\EE^{1/4}\leq\nu\QQ\\
        \Sob{w}{L^\infty}&\leq C\PP^{1/4}\EE^{1/4}\leq\nu\QQ\\
        \Sob{\nabla w}{L^4}&\leq C\PP^{1/4}\ZZ^{1/4}\leq\nu\QQ,
        \end{split}
    \end{align}
and similarly, upon applying \eqref{def:rad:H1} in addition, we have
    \begin{align}\label{est:B:DB}
        \begin{split}
            |B(w,w)|&\leq C\PP^{1/4}\ZZ^{1/2}\EE^{1/4}\leq C\PP\leq\nu^2\QQ\\
            |(DBu)w|&\leq C\nu G\left((\s_1^{1/2}+G)^{1/2}\ZZ^{1/2}+\PP^{1/4}\EE^{1/4}\right)\leq C\nu(\s_1^{1/2}+G)^{3/2}\PP^{1/2}\leq\nu^2\QQ.
        \end{split}
    \end{align}
We are now ready to prove \cref{lem:power}.

\begin{proof}[Proof of \cref{lem:power}]
Observe that $\FF$ from \eqref{def:F} may initially be estimated with Young's inequality and the Cauchy-Schwarz inequality to obtain
    \begin{align}\label{est:FF:partial}
        |\FF+4\mu\EE_N\JJd_1+4\nu\ZZ_N\JJd_1-2\JJ_1\JJd_1|&\leq C\mu\JJ_1^2+C\nu^3\left(\frac{\nu}{\mu}\right)\PP^2+C\nu\left(\frac{\nu}{\mu}\right)\JJ_2^2+\frac{1}{100}\DD\notag\\
        &\leq C\mu\JJ_1^2+C\nu\left(\frac{\nu}{\mu}\right)\JJ_2^2+\nu^7\tQQ%\lam_3\QQ
        +\frac{1}{100}\DD.
    \end{align}
This leaves the remainder of the terms $\FF$ to be
    \begin{align}\label{est:FF:remainder}
        -4\mu\EE_N\JJd_1-4\nu\ZZ_N\JJd_1+2\JJ_1\JJd_1.
    \end{align}
Thus, to complete the estimate of \eqref{def:F}, we must treat $\JJ_1$ and $\JJ_2$, defined in \eqref{def:J}, which will complete the estimate \eqref{est:FF:partial} (Step 1). We will then estimate the remainder of $\FF$ in \eqref{est:FF:remainder}. From \eqref{def:F}, we see that this will require estimates for $\JJd_1$ (Step 2). Once such $\JJd_1$ have been established, we obtain a complete estimate for $\FF$ (Step 3). To these ends, in what follows, we will only track terms whose pre-factors depend increasingly on $\mu$ and hide all others within $\QQ$.

\subsubsection*{Step 1: Estimates $\JJ_1,\JJ_2$}
For $\JJ_1$, we make use of \eqref{eq:ortho:B} to see that
    \begin{align}\notag
        \lb(DB_Nu)w,w_N\rb=\lb B(u,w),w_N\rb+\lb B(w,u),w_N\rb=-\lb B(u,w_N),Q_Nw\rb+\lb B(w,u),w_N\rb,
    \end{align}
where $Q_N=I-P_N$. Then by H\"older's inequality, \eqref{est:diff:interpolation}, the Poincar\'e inequality, inverse Poincar\'e inequality, and \eqref{def:rad:H1}, we estimate $\JJ_1$ as
    \begin{align}
        |\JJ_1|&\leq |\lb B(w,w_N),Q_Nw\rb|+|\lb(DB_Nu)w,w_N\rb|+|\De\nu||\lb A\tu_N,w_N\rb|\notag\\
        &\leq \Sob{w}{L^4}\lVert w_N\rVert\Sob{Q_Nw}{L^4}+\Sob{u}{L^4}\Sob{\nabla w_N}{L^4}|Q_Nw|+\Sob{w}{L^4}\lVert u\rVert\Sob{w_N}{L^4}+|\De\nu|\lVert\tu_N\rVert\lVert w_N\rVert\notag\\
        &\leq \frac{C}{N^{3/2}}\PP\ZZ_N^{1/2}+\frac{C}{N^2}\nu G\PP^{1/2}\PP_N^{1/4}\ZZ_N^{1/4}+C\nu G\PP^{1/2}\ZZ_N^{1/4}\EE_N^{1/4}+\nu^2\left(\frac{|\De\nu|}{\nu}\right)\left(\frac{\lVert\tu\rVert}{\nu}\right)\ZZ_N^{1/2}\notag\\
        &\leq \nu^3\QQ.\label{est:JJ1}
    \end{align}
    
Next, we treat $\JJ_2$. Making use of \eqref{eq:ortho:B} once again, we see that
    \begin{align}
        &\lb B_N(w,w),Aw_N\rb=\lb w^j\bdy_jw^k,\bdy_\ell^2w_N^k\rb=-\lb \bdy_\ell w^j\bdy_jw^k,\bdy_\ell w_N^k\rb+\lb w^j\bdy_j\bdy_\ell w_N^k,\bdy_\ell Q_Nw^k\rb\notag.
    \end{align}
Similarly
    \begin{align}
        \lb(DB_Nu)w,Aw_N\rb&=\lb B(u,w),Aw_N\rb+\lb B(w,u),Aw_N\rb=\lb u^j\bdy_jw^k,\bdy_\ell^2w_N^k\rb+\lb w^j\bdy_ju^k,\bdy_\ell^2w_N^k\rb\notag\\
        &=-\lb \bdy_\ell u^j\bdy_jw^k,\bdy_\ell w_N^k\rb-\lb \bdy_\ell w^j\bdy_ju^k,\bdy_\ell w_N^k\rb.\notag
    \end{align}
As in $\JJ_1$, we then estimate
    \begin{align}
        |\JJ_2|&\leq |\lb B_N(w,w),Aw_N\rb|+|\lb (DB_Nu)w,Aw_N\rb| + |\De\nu||\lb A\tu_N, Aw_N\rb|\notag\\
        &\leq\Sob{\nabla w}{L^4}^2\lVert w_N\rVert+\Sob{w}{L^4}|Aw_N|\Sob{\nabla Q_Nw}{L^4}+2\lVert u\rVert\Sob{\nabla w}{L^4}\Sob{\nabla w_N}{L^4}+|\De\nu||A\tu_N||Aw_N|\notag\\
        &\leq C\PP^{3/2}+C\nu G\PP+\nu^2\left(\frac{|\De\nu|}{\nu}\right)\left(\frac{|A\tu_N|}{\nu}\right)\PP^{1/2}\notag\\
        &\leq \nu^3\QQ\label{est:JJ2}
    \end{align}
    
\subsubsection*{Step 2: Estimates for $\JJd_1$}
First, let us expand $\JJd_1$ so that
    \begin{align}\label{eq:JJ1d:expand}
        \JJd_1&=\lb B(\frac{dw}{dt},w_N),Q_Nw\rb+\lb B(w,\frac{dw_N}{dt}),Q_Nw\rb+\lb B(w,w_N),\frac{dQ_Nw}{dt}\rb\notag\\
        &\quad-\lb (DB_N\frac{du}{dt})w,w_N\rb-\lb(DB_Nu)\frac{dw}{dt},w_N\rb-\lb (DB_Nu)w,\frac{dw_N}{dt}\rb\notag\\
        &\quad-(\De\nu)\lb \frac{d\tu_N}{dt},Aw_N\rb-(\De\nu)\lb A\tu_N,\frac{dw_N}{dt}\rb\notag\\
        &=\II_1+\II_2+\II_3+\KK_1+\KK_2+\KK_3+\LL_1+\LL_2.
    \end{align}
Observe that
    \begin{align}\label{eq:JJ1d:evols}
    \begin{split}
     \frac{dw}{dt}+\nu Aw&=-B(w,w)-(DBu)w-(\De\nu) A\tu-\mu w_N\\
     \frac{dw_N}{dt}+\nu Aw_N&=-B_N(w,w)-(DB_Nu)w-(\De\nu) A\tu_N-\mu w_N\\
     \frac{dQ_Nw}{dt}+\nu AQ_Nw&=-Q_NB(w,w)-(Q_NDBu)w-(\De\nu) AQ_N\tu\\
    \frac{d\tu_N}{dt}+\tnu A\tu_N&=-B_N(\tu,\tu)+f_N-\mu w_N.
    \end{split}
    \end{align}
We may then expand $\II_1$--$\II_3$ as
    \begin{align}
       \II_1=&-\nu\lb B(Aw,w_N),Q_Nw\rb-\lb B(B(w,w),w_N),Q_Nw\rb-\lb B((DBu)w,w_N),Q_Nw\rb\notag\\
        &-(\De\nu)\lb B(A\tu,w_N),Q_Nw\rb-\mu\lb B(w_N,w_N),Q_Nw\rb\notag\\
        \II_2=&-\nu\lb B(w,Aw_N),Q_Nw\rb-\lb B(w,B_N(w,w)),Q_Nw\rb-\lb B(w,(DB_Nu)w),Q_Nw\rb\notag\\
        &-(\De\nu)\lb B(w,A\tu_N),Q_Nw\rb-\mu\lb B(w,w_N),Q_Nw\rb\notag\\
        \II_3=&-\nu\lb B(w,w_N),AQ_Nw\rb-\lb B(w,w_N),Q_NB(w,w)\rb-\lb B(w,w_N),(Q_NDBu)w\rb\notag\\
        &-(\De\nu)\lb B(w,w_N),AQ_N\tu\rb.\notag
    \end{align}
Similarly, we expand $\KK_1$--$\KK_3$ as
    \begin{align}
        \KK_1=&\nu\lb(DB_N( Au)w,w_N\rb+\lb(DB_N\left(B(u,u)\right))w,w_N\rb-\lb(DB_Nf)w,w_N\rb\notag\\
        \KK_2=&\nu\lb(DB_Nu)Aw,w_N\rb+\lb(DB_Nu)B(w,w),w_N\rb+\lb(DB_Nu)^2w,w_N\rb\notag\\
        &+(\De\nu)\lb(DB_Nu)A\tu_N,w_N\rb+\mu\lb(DB_Nu)w_N,w_N\rb\notag\\
        \KK_3=&\nu\lb(DB_Nu)w,Aw_N\rb+\lb(DB_Nu)w,B_N(w,w)\rb+|(DB_Nu)w|^2\notag\\
        &+(\De\nu)\lb(DB_Nu)w,A\tu_N\rb+\mu\lb(DB_Nu)w,w_N\rb.\notag
    \end{align}
Lastly, we expand $\LL_1,\LL_2$ as
    \begin{align}
        \LL_1=&\tnu(\De\nu)\lb A\tu_N,Aw_N\rb+(\De\nu)\lb B_N(\tu,\tu),Aw_N\rb-(\De\nu)\lb f_N,Aw_N\rb+\mu(\De\nu)\lVert w_N\rVert^2\notag\\
        \LL_2=&\nu(\De\nu)\lb A\tu_N,Aw_N\rb+(\De\nu)\lb A\tu_N,B_N(w,w)\rb+(\De\nu)\lb A\tu_N,(DB_Nu)w\rb\notag\\
        &+(\De\nu)^2|A\tu_N|^2+\mu(\De\nu)\lb A\tu_N,w_N\rb.\notag
    \end{align}
We may thus extract the sign-definite terms to rewrite $\JJd_1$ as
    \begin{align}\label{eq:JJ1d:expand:equiv}
        -\JJd_1&+|(DB_Nu)w|^2+(\De\nu)^2|A\tu_N|^2=-\II_1-\II_2-\II_3-\KK_1-\KK_2-\tKK_3-\LL_1-\tLL_2,
    \end{align}
and in the particular case, when $\De\nu\geq0$, we have
    \begin{align}\label{eq:JJ1d:expand:equiv:mod}
        -\JJd_1+|(DB_Nu)w|^2+(\De\nu)^2|A\tu_N|^2&+\mu(\De\nu)\lVert w_N\rVert^2\notag\\
        &=-\II_1-\II_2-\II_3-\KK_1-\KK_2-\tKK_3-\tLL_1-\tLL_2,
    \end{align}
where
    \begin{align}\label{def:tKtL}
        \begin{split}
        \tKK_3&:=\KK_3-|(DB_Nu)w|^2,\\
        \tLL_1&:=\LL_1-\mu(\De\nu)\lVert w_N\rVert^2\\
        \tLL_2&:=\LL_2-(\De\nu)^2|A\tu_N|^2.
        \end{split}
    \end{align}
We will obtain estimates for the right-hand side of \eqref{eq:JJ1d:expand:equiv} and \eqref{eq:JJ1d:expand:equiv:mod}.

We first estimate $\II_1$--$\II_3$ from \eqref{eq:JJ1d:expand}. We estimate $\II_1$ with \eqref{est:diff:interpolation}, \eqref{est:B:DB}, the Poincar\'e inequality, the inverse Poincar\'e inequality, and Young's inequality to obtain
    \begin{align}
        |\II_1|&\leq\left(\nu|Aw|+|B(w,w)|\right)\Sob{\nabla w_N}{L^4}\Sob{Q_Nw}{L^4}+|(DBu)w|\Sob{\nabla w_N}{L^4}\Sob{Q_Nw}{L^4}\notag\\
         &\quad+|\De\nu||A\tu_N|\Sob{\nabla w_N}{L^4}\Sob{Q_Nw}{L^4}+\mu\Sob{w_N}{L^\infty}\lVert w_N\rVert\lVert Q_Nw\rVert\notag\\
         &\leq\frac{C}{N^{3/2}}\left(\nu\PP+\PP^{3/2}\right)\PP_N^{1/4}\ZZ_N^{1/4}+\frac{C}{N^{3/2}}\nu(\s_1^{1/2}+G)^{3/2}\PP\PP_N^{1/4}\ZZ_N^{1/4}\notag\\%+\frac{C}{N^3}\PP\PP_N^{1/2}\ZZ_N^{1/2}\notag\\
         &\quad+\frac{1}{100}(\De\nu)^2|A\tu_N|^2+\frac{C}{N^3}\PP\PP_N^{1/2}\ZZ_N^{1/2}+\frac{C}{N}\mu\PP^{1/2}\PP_N^{1/4}\ZZ_N^{1/2}\EE_N^{1/4}\notag\\
         &\leq \nu^4\QQ+\frac{C}{N}\mu\PP^{1/2}\PP_N^{1/4}\ZZ_N^{1/2}\EE_N^{1/4}+\frac{1}{100}(\De\nu)^2|A\tu_N|^2\notag.
    \end{align}
We treat $\II_2$ and $\II_3$ similarly. For $\II_2$, we have
    \begin{align}
        |\II_2|&\leq\left(\nu|Aw_N|+|B_N(w,w)|\right)\Sob{w}{L^\infty}\lVert Q_Nw\rVert+|(DB_Nu)w|\Sob{w}{L^\infty}\lVert Q_Nw\rVert\notag\\
        &\quad+|\De\nu||A\tu_N|\Sob{w}{L^\infty}\lVert Q_Nw\rVert+\mu\Sob{w}{L^\infty}\lVert w_N\rVert|Q_Nw|\notag\\
        &\leq \frac{C}{N}\left(\nu\PP_N^{1/2}+\PP^{1/4}\ZZ^{1/2}\EE^{1/4}\right)\PP^{3/4}\EE^{1/4}
        +\frac{1}{100}|(DB_Nu)w|^2+\frac{C}{N^2}\PP^{3/2}\EE^{1/2}\notag\\%+\frac{C}{N^2}\nu^4(\s_1^{1/2}+G)^{3/4}\PP^{5/4}\EE^{1/4}\notag\\
        &\quad+\frac{1}{100}(\De\nu)^2|A\tu_N|^2+\frac{C}{N^2}\mu\PP^{3/4}\EE^{1/4}\ZZ_N^{1/2}\notag\\
        &\leq \nu^4\QQ+\frac{C}{N^2}\mu\PP\ZZ_N^{1/2}+\frac{1}{100}|(DB_Nu)w|^2+\frac{1}{100}(\De\nu)^2|A\tu_N|^2\notag.
    \end{align}
For $\II_3$, we have
    \begin{align}
        |\II_3|&\leq\left(\nu|AQ_Nw|+|Q_NB(w,w)|\right)\Sob{w}{L^\infty}\lVert w_N\rVert+|(Q_NDBu)w|\Sob{w}{L^\infty}\lVert w_N\rVert\notag\\
        &\quad+|\De\nu||AQ_N\tu|\Sob{w}{L^\infty}\lVert w_N\rVert\notag\\
        &\leq C\left(\nu\PP_N^{1/2}+\frac{1}N\left(\Sob{\nabla w}{L^4}^2+\Sob{w}{L^\infty}|Aw|\right)\right)\PP^{1/2}\ZZ_N^{1/2}\notag\\
        &\quad+\frac{C}N\left(\Sob{\nabla u}{L^4}\Sob{\nabla w}{L^4}+\Sob{u}{L^\infty}|Aw|\right)\PP^{1/2}\ZZ_N^{1/2}%\notag\\
        %&\quad
        +\nu^2\frac{C}{N}\left(\frac{|\De\nu|}{\nu}\right)\left(\frac{|A^{3/2}\tu|}{\nu}\right)\PP^{1/4}\EE^{1/4}\ZZ_N^{1/2}\notag\\%+C\PP^{1/2}\EE^{1/2}\ZZ_N+\frac{1}{100}(\De\nu)^2\notag\\
        &\leq \nu^4\QQ\notag.%\nu^4\tQQ+\frac{\nu^4}{N}\left(\frac{|\De\nu|}{\nu}\right)\QQ\notag.%\frac{1}{100}(\De\nu)^2|A\tu_N|^2\notag.
    \end{align}
We combine the estimates above, invoke \eqref{cond:mu:N:nu}, \eqref{def:rad:H2:ng}, and the fact that $N\geq1$ to arrive at
    \begin{align}\label{est:summary:II}
        \sum_{j=1}^3|\II_j|&\leq \nu^4\QQ+
        C\nu^{1/2}\mu^{1/2}\PP^{3/2}+\frac{1}{100}|(DB_Nu)w|^2+\frac{2}{100}(\De\nu)^2|A\tu_N|^2.%+\frac{\nu^4}{N}\left(\frac{|\De\nu|}{\nu}\right)\QQ.
    \end{align}
%Note that Poincar\'e inequality was not applied when pre-factors were independent of $N^{-1}$.

Next, we treat the terms $\KK_1, \KK_2$ and $\tKK_3$ from \eqref{eq:JJ1d:expand} and \eqref{def:tKtL}. They are all estimated with H\"older's inequality, the Poincar\'e inequality, \eqref{est:diff:interpolation}, \eqref{est:B:DB}, and integrating by parts, as needed. We obtain
    \begin{align}
        &|\KK_1|\notag\\
        &\leq \left(\nu|Au|\Sob{w_N}{L^\infty}+|B(u,u)|\Sob{w_N}{L^\infty}+\Sob{f}{L^\infty}|w_N|\right)\lVert w\rVert\notag\\
        &\quad+\left(\nu \lVert w_N\rVert|Au|+\lVert w_N\rVert|B(u,u)|+\lVert f\rVert|w_N|\right)\Sob{w}{L^\infty}\notag\\
        &\leq C\nu^2\left(\left(\frac{|Au|}{\nu}\right)\PP_N^{1/4}\EE_N^{1/4}+\left(\frac{|Au|^{1/2}\lVert u\rVert|u|^{1/2}}{\nu^2}\right)\PP_N^{1/4}\EE_N^{1/4}+\left(\frac{|Af|^{1/2}|f|^{1/2}}{\nu^2}\right)\EE_N^{1/2}\right)\ZZ^{1/2}
        \notag\\
        &\quad+C\nu^2\left(\left(\frac{|Au|}{\nu}\right)\ZZ_N^{1/2}+\left(\frac{|Au|^{1/2}|\lVert u\rVert|u|^{1/2}}{\nu^2}\right)\ZZ_N^{1/2}+\left(\frac{\lVert f\rVert}{\nu^2}\right)\EE_N^{1/2}\right)\PP^{1/4}\EE^{1/4}\notag\\
        &\leq \nu^4\QQ\notag.
    \end{align} 
We estimate each term of $\lb(DB_Nu)\frac{dw}{dt},w_N\rb$, integrating by parts as needed, as
    \begin{align}
        \nu|\lb(DB_Nu)Aw,w_N\rb|&\leq\nu\left(\Sob{u}{L^\infty}\lVert w_N\rVert+\lVert u\rVert\Sob{w_N}{L^\infty}\right)|Aw|\notag\\
            &\leq C\nu^2\left(\left(\frac{|Au|^{1/2}|u|^{1/2}}{\nu}\right)\ZZ_N^{1/2}+\left(\frac{\lVert u\rVert}{\nu}\right)\PP_N^{1/4}\EE_N^{1/4}\right)\PP^{1/2}\notag\\
            &\leq \nu^4\QQ\notag\\
        |\lb (DB_Nu)B(w,w),w_N\rb|&\leq \left(\Sob{u}{L^\infty}\lVert w_N\rVert+\lVert u\rVert\Sob{w_N}{L^\infty}\right)|B(w,w)|\notag\\
        &\leq C\nu\left(\left(\frac{|Au|^{1/2}|u|^{1/2}}{\nu}\right)\ZZ_N^{1/2}+\left(\frac{\lVert u\rVert}{\nu}\right)\PP_N^{1/4}\EE_N^{1/4}\right)\notag\\
        &\quad\times\PP^{1/4}\ZZ^{1/2}\EE^{1/4}\notag\\
        &\leq \nu^4\QQ\notag\\
        |\lb(DB_Nu)(DB_Nu)w,w_N\rb|&\leq |\lb B(u,(DB_Nu)w),w_N\rb|+|\lb B((DB_Nu)w,u),w_N\rb|\notag\\
        &\leq \Sob{u}{L^\infty}\left(\lVert w_N\rVert+|w_N|\right)|(DB_Nu)w|\notag\\
        &\leq \nu^4\QQ+\frac{1}{100}|(DB_Nu)w|^2\notag\\
        |\De\nu||\lb(DB_Nu)A\tu_N,w_N\rb|&\leq |\De\nu|\left(|\lb B_N(u,w_N),A\tu_N\rb|+|\lb B_N(A\tu_N,u),w_N\rb|\right)\notag\\
        &\leq|\De\nu|\left(\Sob{u}{L^\infty}\lVert w_N\rVert+\lVert u\rVert\Sob{w_N}{L^\infty}\right)|A\tu_N|\notag\\
        &\leq C\nu|\De\nu|\left[\left(\frac{|Au|^{1/2}|u|^{1/2}}{\nu}\right)\ZZ_N^{1/2}+\left(\frac{\lVert u\rVert}{\nu}\right)\PP_N^{1/4}\EE_N^{1/4}\right]|A\tu_N|\notag\\
        &\leq \nu^4\QQ+\frac{1}{100}(\De\nu)^2|A\tu_N|^2\notag\\
        \mu\lb(DB_Nu)w_N,w_N\rb&\leq \mu\lVert u\rVert\Sob{w_N}{L^4}^2\leq C\mu\nu G\ZZ_N^{1/2}\EE_N^{1/2}\notag.
    \end{align}
In particular, collecting the estimates, then making a final application of Young's inequality and the Poincar\'e inequality yields
    \begin{align}
       |\KK_2|\leq \nu^4\QQ+C\mu\nu G\ZZ_N^{1/2}\EE_N^{1/2}+\frac{1}{100}|(DB_Nu)w|^2+\frac{1}{100}(\De\nu)^2|A\tu_N|^2.\notag
    \end{align}
Now we treat $\tKK_3$. Similar to the previous two terms, we obtain
    \begin{align}
    \nu|\lb(DB_Nu)w,Aw_N\rb|&\leq\nu|(DB_Nu)w||Aw_N|\leq \nu^4\QQ+\frac{1}{100}|(DB_Nu)w|^2\notag\\
    |\lb(DB_Nu)w,B_N(w,w)\rb|&\leq|(DB_Nu)w||B_N(w,w)|\leq \nu^4\QQ+\frac{1}{100}|(DB_Nu)w|^2\notag\\
    |\De\nu||\lb(DB_Nu)w,A\tu_N\rb|&\leq|\De\nu||(DB_Nu)w||A\tu_N|\leq \nu^4\QQ+\frac{1}{100}(\De\nu)^2|A\tu_N|^2\notag\\
    \mu|\lb(DB_Nu)w,w_N\rb|&\leq \mu|(DB_Nu)w||w_N|\leq C\mu\nu(\s_1^{1/2}+G)^{3/2}\PP^{1/2}\EE_N^{1/2}.\notag
    \end{align}
Upon collecting these estimates, then applying Young's inequality and the Poincar\'e inequality, we obtain
    \begin{align}
        |\tKK_3|\leq \nu^4\QQ+C\mu\nu(\s_1^{1/2}+G)^{3/2}\PP^{1/2}\EE_N^{1/2}+\frac{2}{100}|(DB_Nu)w|^2+\frac{1}{100}(\De\nu)^2|A\tu_n|^2.\notag
    \end{align}
Combining the estimates for $\KK_1,\KK_2,\tKK_3$ and a final application of Young's inequality yields
    \begin{align}\label{est:summary:KK}
        &|\KK_1|+|\KK_2|+|\tKK_3|\notag\\
        &\leq\nu^4\QQ+C\mu\nu G\ZZ_N^{1/2}\EE_N^{1/2}+C\mu\nu(\s_1^{1/2}+G)^{3/2}\PP^{1/2}\EE_N^{1/2}+\frac{3}{100}|(DB_Nu)w|^2+\frac{2}{100}(\De\nu)^2|A\tu_N|^2.
    \end{align}

Finally, we estimate $\LL_1$ and $\tLL_1,\tLL_2$ from \eqref{def:tKtL}. Observe that
    \begin{align}
        |\LL_1|&\leq \tnu|\De\nu||A\tu_N||Aw_N|+|\De\nu||B_N(\tu,\tu)||Aw_N|+|\De\nu|\lVert f_N\rVert\lVert w_N\rVert+2\mu|\De\nu|\ZZ_N\notag\\
        &\leq C\nu^2\left(\frac{\tnu}{\nu}\right)^2\PP_N+\frac{1}{100}(\De\nu)^2|A\tu_N|^2+ |\De\nu|\left(C\Sob{\tu}{L^\infty}\lVert\tu\rVert\PP_N^{1/2}+C\s_1\nu^2 G\ZZ_N^{1/2}+2\mu\ZZ_N\right)\notag\\
        &\leq \nu^4\QQ+2\mu\nu\left(\frac{|\De\nu|}{\nu}\right)\ZZ_N+\frac{1}{100}(\De\nu)^2|A\tu_N|^2\notag.
    \end{align}
On the other hand
    \begin{align}
        |\tLL_1|&\leq \nu^4\QQ+\frac{1}{100}(\De\nu)^2|A\tu_N|^2\notag.
    \end{align}
Similarly
    \begin{align}
        |\tLL_2|&\leq |\De\nu|\left(\nu|Aw_N|+|B_N(w,w)|+|(DB_Nu)w|+{\mu|w_N|}\right)|A\tu_N|\notag\\
        &\leq \nu^4\QQ+C\mu\nu^2\left(\frac{|\De\nu|}{\nu}\right)\left(\frac{|A\tu_N|}{\nu}\right)\EE_N^{1/2}+\frac{1}{100}(\De\nu)^2|A\tu_N|^2\notag.
    \end{align}
Combining and $\LL_1,\tLL_2$ and $\tLL_1, \tLL_2$, we obtain   
    \begin{align}
        |\LL_1|+|\tLL_2|&\leq \nu^4\QQ+2\mu\nu\left(\frac{|\De\nu|}{\nu}\right)\ZZ_N+C\mu\nu^2\left(\frac{|\De\nu|}{\nu}\right)\left(\frac{|A\tu_N|}{\nu}\right)\EE_N^{1/2}+\frac{2}{100}(\De\nu)^2|A\tu_N|^2\label{est:summary:LL}\\
           |\tLL_1|+|\tLL_2|&\leq \nu^4\QQ+C\mu\nu^2\left(\frac{|\De\nu|}{\nu}\right)\left(\frac{|A\tu_N|}{\nu}\right)\EE_N^{1/2}+\frac{2}{100}(\De\nu)^2|A\tu_N|^2.\label{est:summary:LtL}
    \end{align}

Finally, we return to \eqref{eq:JJ1d:expand:equiv}. We apply the estimates \eqref{est:summary:II}, \eqref{est:summary:KK}, \eqref{est:summary:LL}, in addition the assumption \eqref{cond:mu:N:nu} and $\mu\geq1$, so that
    \begin{align}\label{est:JJd1:neg}
    -\JJd_1&\leq \nu^4\QQ+C\nu^{1/2}\mu^{1/2}\PP^{3/2}+C\mu\nu G\ZZ_N^{1/2}\EE_N^{1/2}+C\mu\nu(\s_1^{1/2}+G)^{3/2}\PP^{1/2}\EE_N^{1/2}\notag\\
    &\quad +2\mu\nu\left(\frac{|\De\nu|}{\nu}\right)\ZZ_N+C\mu\nu^2\left(\frac{|\De\nu|}{\nu}\right)\left(\frac{|A\tu_N|}{\nu}\right)\EE_N^{1/2}.
    \end{align}
When $\De\nu\geq0$, then we apply \eqref{est:summary:II}, \eqref{est:summary:KK}, and \eqref{est:summary:LtL} in \eqref{eq:JJ1d:expand:equiv:mod} so that
    \begin{align}\label{est:JJd1:neg:mod}
    -\JJd_1&\leq \nu^4\QQ+C\nu^{1/2}\mu^{1/2}\PP^{3/2}+C\mu\nu G\ZZ_N^{1/2}\EE_N^{1/2}+C\mu\nu(\s_1^{1/2}+G)^{3/2}\PP^{1/2}\EE_N^{1/2}\notag\\
    &\quad +C\mu\nu^2\left(\frac{|\De\nu|}{\nu}\right)\left(\frac{|A\tu_N|}{\nu}\right)\EE_N^{1/2}.
    \end{align}
On the other hand, from \eqref{eq:JJ1d:expand:equiv}, the estimates above, and \eqref{est:B:DB}, we also have
    \begin{align}\label{est:JJd1:pos}
        \left|\JJd_1\right|&\leq \nu^4\QQ+C\nu^{1/2}\mu^{1/2}\PP^{3/2}+C\mu\nu G\ZZ_N^{1/2}\EE_N^{1/2}+C\mu\nu(\s_1^{1/2}+G)^{3/2}\PP^{1/2}\EE_N^{1/2}\notag\\
    &\quad +2\mu\nu\left(\frac{|\De\nu|}{\nu}\right)\ZZ_N+C\mu\nu^2\left(\frac{|\De\nu|}{\nu}\right)\left(\frac{|A\tu_N|}{\nu}\right)\EE_N^{1/2}+{C\nu^2\left(\frac{|\De\nu|}{\nu}\right)^2|A\tu_N|^2}.
    \end{align}

\subsubsection*{Step 3: Completing estimate of $\FF$}

Firstly, from \eqref{est:JJ1}, \eqref{est:JJ2}, and the fact that $\mu\geq1$, it follows that
    \begin{align}\label{est:J:squared}
        C\mu\JJ_1^2+C\nu\left(\frac{\nu}{\mu}\right)\JJ_2^2\leq \mu\nu^6\tQQ,%(\QQ+1)\PP
    \end{align}
where we invoked the notation \eqref{def:tQQ} and $\tQQ$ is to be interpreted as being of the first form indicated there. On the other hand, recalling the definition of $\FF$ in \eqref{def:F}, we see that must also obtain upper bounds for the following terms
    \begin{align}\notag
        -4\mu\EE_N\JJd_1,\quad -{4\nu}\ZZ_N\JJd_1,\quad -2\JJ_1\JJd_1.
    \end{align}
Then we claim that from \eqref{est:JJd1:neg}, after multiple applications of Young's inequality and the Cauchy-Schwarz inequality, that independently of the sign of $\De\nu$, we obtain
    \begin{align}\label{est:JJd1:RHS:a}
        -4\mu\EE_N\JJd_1-4\nu\ZZ_N\JJd_1\leq \nu^7\tQQ+\left[\frac{1}{10}+\frac{7}{12}\left(\frac{|\De\nu|}{\nu}\right)+C\left(\frac{\nu}{\mu}\right)^2G^4\right]\DD.
    \end{align}
where we have also applied the definition of \eqref{def:D}. On the other hand, when $\De\nu\geq0$ is known, then
    \begin{align}\label{est:JJd1:RHS:b}
        -4\mu\EE_N\JJd_1-4\nu\ZZ_N\JJd_1\leq \nu^7\tQQ+\left[\frac{1}{10}+C\left(\frac{\nu}{\mu}\right)^3G^4\right]\DD.
    \end{align}

To see \eqref{est:JJd1:RHS:a}, we multiply the right-hand side of \eqref{est:JJd1:neg} by $\mu\EE_N$, then estimate each of the terms in the resulting product. We obtain
    \begin{align}
        \mu\nu^4\QQ\EE_N&\leq C\nu^7\left(\frac{\nu}{\mu}\right)\QQ^2+\frac{1}{1000}\mu^3\EE_N^2\leq \nu^7\tQQ+\frac{1}{1000}\DD\notag\\
        C\nu^{1/2}\mu^{3/2}\PP^{3/2}\EE_N&\leq C\nu\PP^3+\frac{1}{1000}\mu^3\EE_N\leq \nu^7\tQQ+\frac{1}{1000}\DD\notag\\
        C\mu^2\nu G\ZZ_N^{1/2}\EE_N^{3/2}&\leq C\frac{\nu^4}{\mu}G^4\ZZ_N^2+\frac{1}{1000}\mu^3\EE_N^2\leq\left(C\nu^2\left(\frac{\nu}{\mu}\right)^2G^4+\frac{1}{1000}\right)\DD\notag\\
        C\mu^2\nu(\s_1^{1/2}+G)^{3/2}\PP^{1/2}\EE_N^{3/2}&\leq C\frac{\nu^8}{\mu}(\s_1^{1/2}+G)^6\PP^2+\frac{1}{1000}\mu^3\EE_N^2\leq \nu^7\tQQ+\frac{1}{1000}\DD\notag\\
        8\mu^2\nu\left(\frac{|\De\nu|}{\nu}\right)\ZZ_N\EE_N&\leq \frac{1}4\left(\frac{|\De\nu|}{\nu}\right)\DD\label{est:Dnu:1}\\
        C\mu^2\nu^2\left(\frac{|\De\nu|}{\nu}\right)\left(\frac{|A\tu_N|}{\nu}\right)\EE_N^{3/2}&\leq C\frac{\nu^8}{\mu}\left(\frac{|\De\nu|}{\nu}\right)^4\left(\frac{|A\tu_N|}{\nu}\right)^4+\frac{1}{1000}\mu^3\EE_N^2\leq \nu^7\tQQ+\frac{1}{1000}\DD\notag.
    \end{align}
Repeating the same argument above, except after multiplying \eqref{est:JJd1:neg} by $\nu\ZZ_N$, we obtain
    \begin{align}
        \nu^5\QQ\ZZ_N&\leq C\frac{\nu^8}{\mu}\QQ^2+\frac{1}{1000}\mu\nu^2\ZZ_N^2\leq \nu^7\tQQ+\frac{1}{1000}\DD\notag\\
        C\nu^{3/2}\mu^{1/2}\PP^{3/2}\ZZ_N&\leq C\nu\PP^3+\frac{1}{1000}\mu\nu^2\ZZ_N^2\leq\nu^7\tQQ+\frac{1}{1000}\DD\notag\\
        C\mu\nu^2G\ZZ_N^{3/2}\EE_N^{1/2}&\leq C\mu\nu^2 G^4\EE_N^2+\frac{1}{1000}\mu\nu^2\ZZ_N^2\leq \left(C\left(\frac{\nu}{\mu}\right)^2G^4+\frac{1}{1000}\right)\DD\notag\\
        C\mu\nu^2(\s_1^{1/2}+G)^{3/2}\PP^{1/2}\EE_N^{1/2}\ZZ_N&\leq C\mu\nu^2(\s_1^{1/2}+G)^3\PP\EE_N+\frac{1}{1000}\mu\nu^2\ZZ_N^2\notag\\
        &\leq C\frac{\nu^4}{\mu}(\s_1^{1/2}+G)^6\PP^2+\frac{1}{1000}\mu^3\EE_N^2+\frac{1}{1000}\DD\notag\\
        &\leq \nu^7\tQQ+\frac{2}{1000}\DD \notag\\
        8\mu\nu^2\left(\frac{|\De\nu|}{\nu}\right)\ZZ_N^2&\leq \frac{1}3\left(\frac{|\De\nu|}{\nu}\right)\DD\label{est:Dnu:2}\\
        C\mu\nu^3\left(\frac{|\De\nu|}{\nu}\right)\left(\frac{|A\tu_N|}{\nu}\right)\EE_N^{1/2}\ZZ_N&\leq C\mu\nu^4\left(\frac{|\De\nu|}{\nu}\right)^2\left(\frac{|A\tu_N|}{\nu}\right)^2\EE_N+\frac{1}{1000}\mu\nu^2\ZZ_N^2\notag\\
        &\leq  C\frac{\nu^8}{\mu}\left(\frac{|\De\nu|}{\nu}\right)^4\left(\frac{|A\tu_N|}{\nu}\right)^4+\frac{1}{1000}\mu^3\EE_N^2+\frac{1}{1000}\DD\notag\\
        &\leq \nu^7\tQQ+\frac{2}{1000}\DD.\notag
    \end{align}
Combining these estimates yields \eqref{est:JJd1:RHS:a}. On the other hand, if $\De\nu\geq0$, then we see from \eqref{est:JJd1:pos} that we instead obtain \eqref{est:JJd1:RHS:b} by simply ignoring the estimates \eqref{est:Dnu:1}, \eqref{est:Dnu:2} above.

Lastly, we see from \eqref{eq:JJ1d:expand:equiv} that we are now left to estimate $-2\JJ_1\JJd_1$. From \eqref{est:JJ1} and \eqref{est:JJd1:pos}, we apply the Poincar\'e inequality to deduce
    \begin{align}\label{est:JJ1:JJd1}
        -2\JJ_1\JJd_1\leq 2|\JJ_1||\JJd_1|\leq \nu^7\tQQ+\mu\nu^6\QQ\PP
    \end{align}
where the second term represents all terms from the resulting product that appear with $\mu$-dependent pre-factors. Note that the dependence on $\mu$ to leading order can be made linear due to the fact that $\mu\geq1$.

Combining \eqref{est:FF:partial}, \eqref{est:JJd1:RHS:a}, and \eqref{est:JJ1:JJd1}, we have
    \begin{align}\label{est:FF:1}
        \FF\leq \nu^7\tQQ+\mu\nu^6\QQ\PP+\left[\frac{2}{10}+\frac{7}{12}\left(\frac{|\De\nu|}{\nu}\right)+C\left(\frac{\nu}{\mu}\right)^2G^4\right]\DD
    \end{align}
On the other hand, when $\De\nu\geq0$, we have
    \begin{align}\label{est:FF:2}
       \FF\leq \nu^7\tQQ+ \mu\nu^6\QQ\PP+\left[\frac{2}{10}+C\left(\frac{\nu}{\mu}\right)^2G^4\right]\DD.
    \end{align}
    
\subsubsection*{Proof of the main claim}
Finally, now we return to \eqref{eq:diff:balance:Ndot:compact}. Then from \eqref{est:FF:1}, we have
    \begin{align}\label{est:EEd:gen}
     \frac{d}{dt}\EEd_N^2+2\mu\EEd_N^2\leq\nu^7\tQQ+\mu\nu^6\QQ\PP+\left[\frac{2}{10}+\frac{7}{12}\left(\frac{|\De\nu|}{\nu}\right)+C\left(\frac{\nu}{\mu}\right)^2G^4-1\right]\DD.
    \end{align}
If the first condition in \eqref{cond:diff:nu} holds, then in particular $|\De\nu|\nu^{-1}\leq1$. Also by \eqref{cond:mu:power}, we see that $\sqrt{C}\nu\mu^{-1}G^2\leq1/10$. These two facts together ensure that the prefactor of $\DD$ is negative. On the other hand, by \eqref{cond:mu:N:nu:tnu}, we see that we may apply \cref{thm:ng:bounds} and \cref{thm:sensitivity} in conjunction with \eqref{cond:nu:error} and \eqref{cond:mu:power} to identify a time $\tau_1'\geq\tau_0$ such that
    \begin{align}\label{est:P:longtime}
        \PP(t)\leq \nu^2\left(\frac{\nu}{\mu}\right)M^2K_2^2\leq\nu^2M^2K_2^2,
    \end{align}
for all $t\geq\tau_1'$. Since $\QQ$ is monotonically increasing in $\PP$ and $|\De\nu|\nu^{-1}$, it follows from \eqref{est:P:longtime} and \eqref{cond:nu:error} that 
    \begin{align}\notag
    \frac{d}{dt}\EEd_N^2+\mu\EEd_N^2\leq C\nu^7\left(\frac{|\De\nu|}{\nu}\right)^2K_2^2,
    \end{align}
for some $C$, for all $t\geq\tau_1'$. Note that we have also bounded all the coefficients of $\QQ$ that depend on $|\tu|, \lVert\tu\rVert, |A\tu|$ using \eqref{def:rad:H1:ng}, \eqref{def:rad:H2:ng}, and \eqref{def:rad:Hk:ng} for $k=3$, from \cref{thm:ng:bounds}. By Gronwall's inequality, we then deduce
    \begin{align}\notag
    \EEd_N^2(t)\leq \EEd_N^2(\tau_1')e^{-\mu(t-\tau_1')}+ C\nu^6\left(\frac{\nu}{\mu}\right)\left(\frac{|\De\nu|}{\nu}\right)^2,
    \end{align}
for all $t\geq\tau_1'$, for some constant $C>0$ depending on $M, K_2$. Hence, there exists $\tau_1>0$ such that
    \begin{align}\label{est:EEd:final}
    \EEd_N^2(t)\leq C\nu^6\left(\frac{\nu}{\mu}\right)\left(\frac{|\De\nu|}{\nu}\right)^2,
    \end{align}
for all $t\geq\tau_1$, as desired.

On the other hand, if $\De\nu\geq0$, i.e., the second condition in \eqref{cond:diff:nu}, then we instead apply \eqref{est:FF:2} to deduce
    \begin{align}\label{est:EEd:pos}
     \frac{d}{dt}\EEd_N^2+2\mu\EEd_N^2\leq\nu^7\tQQ+\mu\nu^6\QQ\PP+\left[\frac{2}{10}+C\left(\frac{\nu}{\mu}\right)^2G^4-1\right]\DD.
    \end{align}
Then we may argue as before, but bypass having to invoke the first condition in \eqref{cond:diff:nu}, to ultimately obtain \eqref{est:EEd:final}, for all $t\geq\tau_1$, with $\tau_1\geq\tau_0$ sufficiently large. This completes the proof of \cref{lem:power}.
\end{proof}

\subsection{Proofs of Theorems \ref{thm:converge} and \ref{thm:converge:second}}\label{sect:convergence:proof}

It will be useful to have the following elementary inequalities before proceeding to the proofs. They can be obtained by apply H\"older's inequality, \eqref{est:diff:interpolation}, \eqref{est:B:DB}, the Poincar\'e inequality, and the inverse Poincar\'e inequality. Indeed, we have
    \begin{align}
        |\lb B(w,w_N),Q_Nw\rb|&\leq \frac{C}{N^{3/2}}\ZZ^{1/4}\EE^{1/4}\ZZ_N^{1/2}\PP^{1/2}\leq\frac{C}{N^{3/2}}\PP^{3/2},\label{est:nlt:1}\\
        |\lb (DB_Nu)w,w_N\rb|&\leq C\nu(\s_1^{1/2}+G)^{3/2}\PP.\label{est:nlt:2}
    \end{align}
With \cref{lem:power} also in hand, we are now ready to prove \cref{thm:converge} and \cref{thm:converge:second}.

\begin{proof}[Proof of Theorem \ref{thm:converge}]
The proof will proceed by induction. We first establish the base case, $m=0$.

Suppose that $\nu_0>\nu$ and let $\De\nu_0=\nu_0-\nu$. We denote the relative error by
    \begin{align}\label{def:delta}
     \de_0:=\frac{\De\nu_0}{\nu}>0.
    \end{align}
Let $\tG_0$ be defined by \eqref{def:tG}, but with $\tnu=\nu_0$. Observe that
    \begin{align}\label{eq:nu:nu0:equiv}
       \frac{\nu_0}\nu=\de_0+1.
    \end{align}
Thus
    \begin{align}\label{est:G:G0:equiv}
    \left(\de_0+1\right)^2G_0=G\geq G_0
    \end{align}
where
    \begin{align}\label{def:G0}
    G_0:=\frac{|g|}{\nu_0^2}.
    \end{align}
Suppose that $\mu_1\geq\nu_0$. Since $\nu\leq\nu_0$, it follows that
    \begin{align}\label{est:tG0}
    \tG^2_0\leq 2\left(\de_0+1\right)^4G_0^2.
    \end{align}
In particular, \eqref{est:G:G0:equiv} and \eqref{est:tG0} imply that $G, \tG_0$ can be controlled entirely in terms of the relative error, $\de_0$, and a modified Grashof number, $G_0$, that depends only on the initial viscosity, $\nu_0$, and the forcing $|g|$.

Suppose that \eqref{cond:mu:N:tnu}, \eqref{cond:mu:N:nu}, \eqref{cond:mu:ng:H2}, and $\eqref{cond:mu:ng:Hk}_3$  all hold with $\tnu=\nu_0$, $\mu=\mu_1$, and $N=N_1$. Upon substituting for all instances of $\nu$ using \eqref{eq:nu:nu0:equiv} in \eqref{cond:mu:N:tnu}, \eqref{cond:mu:N:nu}, \eqref{cond:mu:ng:H2}, and $\eqref{cond:mu:ng:Hk}_3$, we see that the conditions on $\mu_1, N_1$ can be described entirely in terms of $\al,\de_0,\s_1,\s_2, G_0$.

Consider the unique global-in-time solutions $u(\cdotp;u_0,\nu)$ and $\tu(\cdotp;\tu_0,\nu_0)$, for some $\al\geq1$, to \eqref{eq:nse:ff} and \eqref{eq:nse:ng:ff} corresponding to initial values $u_0\in \bigcap_{\ell=1}^3B_\ell(R_\ell)$, $\tu_0\in \bigcap_{\ell=1}^3B_\ell(\al R_\ell)$ and viscosities, $\nu,\nu_0$, respectively, and nudging parameter $\mu=\mu_1$. Then also consider the system \eqref{eq:diff:balance:Ndot:I0} on the interval $\til{I}_0=[0,\infty)$. We will choose $t_1>0$ specifically, but as of now it is arbitrary. Let $\nu_1>0$ be given by \eqref{def:update} with $m=0$ and let $\De\nu_1=\nu_1-\nu$.

From \eqref{eq:induction:step} for $\veps=\veps_1$, we have
    \begin{align}\label{est:Dnu:1:prep}
        |\De\nu_1|
        &\leq \frac{1}{\nu_0^2\veps}\left(\left|\EEd_N^{(1)}\right|+2\nu\ZZ_N^{(1)}+\left|\lb B(w^{(1)},w_N^{(1)}),Q_Nw^{(1)}\rb\right|+\left|\lb (DB_Nu^{(1)})w^{(1)},w_N^{(1)}\rb\right|\right).
    \end{align}
%Applying \eqref{est:nlt:1}, we obtain%From H\"older's inequality, \eqref{est:diff:interpolation}, the Poincar\'e inequality, and inverse Poincar\'e inequality, we have
%    \begin{align}\notag
 %       |\lb B(w^{(1)},w_N^{(1)}),Q_Nw^{(1)})\rb|\leq \frac{C}{N^{3/2}}(\ZZ^{(1)})^{1/4}(\EE^{(1)})^{1/4}(\ZZ^{(1)}_N)^{1/2}(\PP^{(1)})^{1/2}\leq \frac{C}{N^{3/2}}(\PP^{(1)})^{3/2}.
    %\end{align}
By \cref{thm:sensitivity}, there exists $\tau_0>0$ such that 
    \begin{align}\label{est:PP:0}
        \PP(t)\leq \nu^2\left(\frac{\nu}{\mu_1}\right)\left(\frac{\De\nu_0}{\nu}\right)\de_0K_2^2,
    \end{align}
for all $t\geq\tau_0$, where $K_2$ is given by \eqref{def:K2}, but in terms of $\tnu=\nu_0$ and $\tG=\tG_0$. Observe that due to \eqref{est:G:G0:equiv}, \eqref{est:tG0}, $K_2$ can be bounded entirely in terms of $\al, \s_1,\s_2,\de_0, G_0$. Hence, for $t_1\geq\tau_0$, we have
    \begin{align}\label{est:ZN:0}
        2\nu\ZZ_N^{(1)}\leq 2\nu^3\left(\frac{\nu}{\mu_1}\right)\left(\frac{\De\nu_0}{\nu}\right)\de_0K_2^2.
    \end{align}
Furthermore, from \eqref{est:nlt:1}, \eqref{est:nlt:2}, \eqref{est:G:G0:equiv} we may also derive
    \begin{align}
     |\lb B(w^{(1)},w_N^{(1)}),Q_Nw^{(1)})\rb|&\leq C\nu^3\left(\frac{\nu}{\mu_1}\right)^{9/4}\left(\frac{\De\nu_0}{\nu}\right)\de_0^2K_2^3,\label{est:B:0}\\
     |\lb (DB_Nu^{(1)})w^{(1)},w_N^{(1)}\rb|
     %&\leq |(DB_Nu^{(1)})w^{(1)}||w_N^{(1)}|\leq C\nu(\s_1^{1/2}+G)^{3/2}\PP^{(1)}\notag\\
        &\leq C\nu^3(\de_0+1)^{3/2}(\s_1^{1/2}+G_0)^{3/2}\left(\frac{\nu}{\mu_1}\right)\left(\frac{\De\nu_0}{\nu}\right)\de_0K_2^2,\label{est:DB:0}
    \end{align}
where we additionally applied \eqref{cond:mu:N:nu}. %Arguing similarly from \eqref{est:nlt:2}, we obtain
 %   \begin{align}\label{est:DB:0}
  %      |\lb (DB_Nu^{(1)})w^{(1)},w_N^{(1)}\rb|&\leq |(DB_Nu^{(1)})w^{(1)}||w_N^{(1)}|\leq C\nu(\s_1^{1/2}+G)^{3/2}\PP^{(1)}\notag\\
   %     &\leq C\nu^3(\s_1^{1/2}+G)^{3/2}\left(\frac{\nu}{\mu}\right)\left(\frac{\De\nu_0}{\nu}\right)\de_0K_2^2.
%    \end{align}
Finally, with the hypothesis of \cref{lem:power} satisfied with $\tnu=\nu_0$, $\mu=\mu_1$, and $N=N_1$, we see that
    \begin{align}\label{est:EEd:0}
        |\EEd_N^{(1)}|\leq \nu^3\left(\frac{\nu}{\mu_1}\right)^{1/2}\left(\frac{\De\nu_0}{\nu}\right)L_0.
    \end{align}
Again, due to \eqref{est:G:G0:equiv}, \eqref{est:tG0}, $L_0$ can be bounded in terms of $\al, \de_0, G_0$. Combining \eqref{est:ZN:0}--\eqref{est:EEd:0} in \eqref{est:Dnu:1:prep}, we obtain
    \begin{align}\label{est:error:1}
        |\De\nu_1|\leq \frac{C}{\veps_1}\left(\frac{1}{\de_0+1}\right)^{1/2}\left(\frac{\nu_0}{\mu_1}\right)^{1/2}\left(L_0+\left(1+\de_0^2K_2+(\de_0+1)^{5/2}(\s_1^{1/2}+G_0)^{3/2}\right)K_2^2\right)\De\nu_0,
    \end{align}
where we applied the fact that $\nu\leq\nu_0\leq\mu_1$. Fix any $\be\in(0,1)$ such that
    \begin{align}\label{cond:beta}
        \frac{\be}{1+\be}<\frac{1}{1+\de_0}.
    \end{align}
Then we choose $\mu_1$ such that
    \begin{align}\label{cond:mu:proof}
    \frac{C}{\veps_1}\left(\frac{1}{\de_0+1}\right)^{1/2}\left(\frac{\nu_0}{\mu_1}\right)^{1/2}\left(L_0+\left(1+\de_0^2K_2+(\de_0+1)^{5/2}(\s_1^{1/2}+G_0)^{3/2}\right)K_2^2\right)\leq\be.
    \end{align}
Upon applying \eqref{cond:mu:proof} in \eqref{est:error:1}, we see that \eqref{cond:beta} ensures that
    \begin{align}\label{eq:positivity}
        \nu_1=(\nu_1-\nu)+\nu\geq(1-\be\de_0)\nu>0.
    \end{align}
This proves that there exists $\be\in(0,1\wedge\de_0^{-1})$  and $t_1>t_0=0$ and $\be\in(0,1)$ such that for $\nu_1$ determined by \eqref{def:update} with $M=1$ one has that $\nu_1>0$ and
    \begin{align}\notag
        |\De\nu_1|\leq\be|\De\nu_0|<\nu.
    \end{align}
Notice that the conditions \eqref{cond:mu:N:tnu}, \eqref{cond:mu:N:nu}, \eqref{cond:mu:ng:H2}, $\eqref{cond:mu:ng:Hk}_3$, and \eqref{cond:mu:proof} can be summarized as
    \begin{align}\notag
        \gam_1\leq\frac{\mu_1}{\nu_0}\leq \gam N_1^2,
    \end{align}
for some $\gam>0$ that depends on $\al,\de_0,\s_1,\s_2,G_0$ and $\gam_1>0$ that depends additionally on $\veps_1^{-1}$. In particular, the base case holds.

Now suppose that for some $M\geq2$, there exist times $0=t_0<\dots, t_{M-1}$ and constants $\gam,\gam_1,\dots, \gam_{M-1}>0$, where $\gam$ depends on $\al,\de_0,\s_1,\s_2, G_0$, and $\gam_m$ depends additionally on $\veps_m^{-1}$, such that if 
    \begin{align}\notag
        \veps_1,\dots,\veps_{M-1}>0,    
    \end{align}
where $\veps_m$ is given by \eqref{def:eps}, and $\mu_1,\dots,\mu_{M-1}>0$, $N_1,\dots, N_{M-1}\geq1$ satisfy
    \begin{align}\label{cond:mu:N:induct}
        \gam_m\leq\frac{\mu_m}{\nu_0}\leq\gam N_m^2,
    \end{align}
for all $m=1,\dots, M-1$, then there exists $\be\in(0,1\wedge\de_0^{-1})$ and $0\leq t_1\leq\dots\leq t_{M-1}$ such that for $\nu_1\dots, \nu_{M-1}$, given by \eqref{def:update}, corresponding to $t_1,\dots, t_{M-1}$, respectively, one has 
    \begin{align}\label{eq:induction:hyp}
        |\De\nu_m|\leq \be|\De\nu_{m-1}|,
    \end{align} 
for all $m=1,\dots, M-1$, and $\nu_1,\dots, \nu_{M-1}>0$. We show that there exists a constant $\gam_M>0$ and $t_M\geq t_{M-1}$ such that for $\nu_M$ defined by \eqref{def:update}, evaluated at time $t=t_M$, one has that $\nu_M>0$ and \eqref{eq:induction:hyp} also holds, whenever $\veps_M>0$ and $\mu_M,N_M$ satisfy \eqref{cond:mu:N:induct} for $m=M$.

Observe that upon iteration of \eqref{eq:induction:hyp}, one has
    \begin{align}\label{eq:iterate}
        |\De\nu_m|\leq \be^{m}\de_0\nu<\nu,
    \end{align}
for all $m=1,\dots, M-1$. In particular
    \begin{align}\label{est:nu:num:equiv}
        0<1-\be\de_0\leq\frac{\nu_m}{\nu}\leq \be\de_0+1.
    \end{align}
Hence, in light of \cref{rmk:equiv} and \eqref{eq:nu:nu0:equiv}, \eqref{est:G:G0:equiv}, we see that any dependence on $\nu$, explicit or implicit, in \eqref{cond:mu:N:nu}, \eqref{cond:mu:ng:H2}, and $\eqref{cond:mu:ng:Hk}_3$, may ultimately be replaced by $\nu_0$ up to a multiplicative factor that depends only on $\be,\de_0$. As a result, all instances of $\tG, G$, where $\tG$ is defined by \eqref{def:tG} with $\tnu=\tnu_m$, is moreover replaced by $G_0$ up to a multiplicative factor that depends only on $\be,\de_0$. Similarly, in \eqref{cond:mu:N:tnu},  $\tnu=\nu_m$ may be replaced by $\nu_0$ up to a multiplicative factor that depends only on $\be,\de_0$.

%In particular, owing to \eqref{est:G:G0:equiv}, we may ultimately replace all instances of $\tG$ with $G_0$, up to a multiplicative constant that depends only on $\be, \de_0$.

Let us now consider $\til{I}_{M-1}=[t_{M-1},\infty)$. We re-initialize the feedback control system over $\til{I}_{M-1}$ with $\tu(t_{M-1}^-)\in B_k(\al R_k)$, where $\tu(t_{M-1}^-)=\tu(t_{M-1};\nu_{M-1})$ and $\al\geq1$ is the same scaling factor used for re-initializing in all previous intervals $\til{I}_0,\til{I}_1,\dots,\til{I}_{M-2}$. Denote the unique global-in-time solutions $u(t;u(t_{M-1}),\nu)$ and $\tu(t;\tu(t_{M-1}^-),\nu_{M})$, for $t\in\til{I}_{M-1}$, to the initial value problems, \eqref{eq:nse:ff} and \eqref{eq:nse:ng:ff} over $\til{I}_{M-1}$, corresponding to initial data $u(t_{M-1};\nu), \tu(t_{M-1}^-;\nu_{M-1})$ and viscosities, $\nu,\nu_{M}$, respectively. From \eqref{def:update}, we have
    \begin{align}
        &|\De\nu_M|\leq \frac{1}{\nu_0^2\veps_M}\left(\left|\EEd_N^{(M)}\right|+2\nu\ZZ_N^{(M)}+\left|\lb B(w^{(M)},w_N^{(M)}),Q_Nw^{(M)}\rb\right|\right.\notag\\
        &\qquad\qquad\qquad\left.+\left|\lb (DB_Nu^{(M)})w^{(M)},w_N^{(M)}\rb\right|\right).\notag
    \end{align}
We estimate exactly in the same way as before and obtain \eqref{est:ZN:0}, \eqref{est:B:0}, \eqref{est:DB:0}, \eqref{est:EEd:0} in terms of $|\De\nu_{M-1}|$ and $\mu_M$ instead, except that we invoke the first condition in \eqref{cond:diff:nu} from \cref{lem:power}, due to \eqref{eq:iterate}, in order to deduce
    \begin{align}\notag
        |\EEd_N(t)|\leq\nu^3\left(\frac{\nu}{\mu}\right)^{1/2}\left(\frac{|\De\nu_{M-1}|}{\nu}\right)L_0,
    \end{align}
for all $t\geq \tau_{M-1}'$, for some $\tau_{M-1}'\geq t_{M-1}$. Note that constants $K_2, L_0$ can again be bounded by quantities that depend only on $\al,\s_1,\s_2,\de_0,G_0$ by applying \eqref{est:G:G0:equiv}, \eqref{def:G0}, \eqref{eq:iterate}, \eqref{est:nu:num:equiv}, and \cref{rmk:equiv}. By choosing to evaluate at $t=t_M>\tau_{M-1}'$, we thus arrive at  
    \begin{align}\notag
     &|\De\nu_M|\\
     &\leq \frac{C}{\veps_M}\left(\frac{1}{\de_0+1}\right)^{1/2}\left(\frac{\nu_0}{\mu_M}\right)^{1/2}\left(L_0+\left(1+\de_0^2K_2+(\de_0+1)^{5/2}(\s_1^{1/2}+G_0)^{3/2}\right)K_2^2\right)|\De\nu_{M-1}|.\notag
    \end{align}
Assuming that $\mu$ satisfies \eqref{cond:mu:proof} with $\veps_1$ replaced by $\veps_M$, where $\be\in(0,1)$ satisfies \eqref{cond:beta}, we may then deduce
    \begin{align}\notag
        |\De\nu_M|\leq \be|\De\nu_{M-1}|.
    \end{align}
Moreover, upon iteration, we have \eqref{eq:iterate}, so that we may again deduce $\nu_M>0$. Lastly, we again observe that the conditions  \eqref{cond:mu:N:nu}, \eqref{cond:mu:ng:H2}, $\eqref{cond:mu:ng:Hk}_3$, and \eqref{cond:mu:proof} with $\veps_1$ replaced by $\veps_M$ can be summarized as
    \begin{align}\notag
        \gam_M\leq\frac{\mu_M}{\nu_0}\leq \gam N_M^2.
    \end{align}
for some $\gam_M, \gam>0$ where $\gam$ depends only on $\al,\de_0,\s_1,\s_2,G_0$, and $\gam_M$ depends additionally on $\veps_M^{-1}$. This completes the proof.
%\qed

\end{proof}

%\subsection{Proof of Theorem \ref{thm:converge:second}}

Next, we prove \cref{thm:converge:second}.

\begin{proof}[Proof of Theorem \ref{thm:converge:second}]
We proceed as in the proof of \cref{thm:converge} above. We consider the base case $m=0$. Let $\bar{\nu}_0>0$ be arbitrary and let 
    \begin{align}\label{def:delta:second}
        \de_0:=\frac{|\De\bar{\nu}_0|}{\nu},\quad \De\bar{\nu}_0:=\bar{\nu}_0-\nu.
    \end{align}
Observe that if $\bar{\nu}_0\geq\nu$, then we have the relations \eqref{eq:nu:nu0:equiv}, \eqref{est:G:G0:equiv}, where $\bar{G}_0$ is given by \eqref{def:G0} with $\bar{\nu}_0$ replacing $\nu_0$. Moreover, if $\mu_1\geq\bar{\nu}_0$, then $\mu\geq\nu$ and we have \eqref{est:tG0}. Otherwise, if $\bar{\nu}_0<\nu$, then we instead have
    \begin{align}\label{eq:nu:nu0:equiv:second}
        \frac{\bar{\nu}_0}{\nu}=1-\de_0,
    \end{align}
and $\de_0<1$. Hence 
    \begin{align}\label{est:G:G0:equiv:second}
        (1-\de_0)^2\bar{G}_0=G\leq \bar{G}_0.
    \end{align}
Moreover, if we suppose that $\mu_1\geq\bar{\nu}_0$, then
    \begin{align}\label{est:tG0:second}
        \tG_0^2\leq 2(1-\de_0)^2G_0^2\leq2(1+\de_0)^4G_0^2.
    \end{align}
In particular, in either case $\bar{\nu}_0\geq\nu$ or $\bar{\nu}_0<\nu$, it  follows that $G,\tG_0$ can be controlled entirely in terms of the relative error, $\de_0$, and a modified Grashof number, $\bar{G}_0$, that depends only on the initial viscosity, $\bar{\nu}_0$, and the forcing $|g|$. 
        
We consider the same initial setup as in the proof of \cref{thm:converge}, with the notation adjusted accordingly. Define $\w_1>0$ by
    \begin{align}\label{def:lam1}
        \w_1^{-1}:=\bar{\nu}_0|J_1|.
    \end{align}
Let $\bar{\nu}_1$ be given by \eqref{def:update:second} with $m=0$ and $J_1:=[t_1',t_1)$. By \eqref{eq:induction:step:second} if $\bar{\veps}_1>0$, it follows that
    \begin{align}
     |\De\bar{\nu}_{1}|\leq& \frac{1}{\bar{\veps}_1\bar{\nu}_0^2}\left(\frac{|\EE_N(t_{1})-\EE_N(t_{1}')|}{|J_{1}|}\right.\notag\\
        &\left.+\frac{1}{|J_{1}|}\int_{J_{1}}2\nu\ZZ_N(s)+|\lb B(w(s),w_N(s)),Q_Nw(s))\rb|+|\lb (DB_Nu(s))w(s),w_N(s)\rb| ds\right).\notag
    \end{align}

From \eqref{est:sensitivity:L2}, the Poincar\'e inequality, \eqref{est:tG0}, \eqref{est:tG0:second}, we see that there exists $\tau_0>0$ such that
    \begin{align}\label{est:EN:0:second}
        \EE_N(t)\leq C\tal_1^2\nu^2\left(\frac{\nu}{\bar{\mu}_1}\right)\left(\frac{|\De\bar{\nu}_0|}{\nu}\right)\de_0(1+\de_0)^4\bar{G}_0^2,
    \end{align}
for all $t\geq \tau_0$, where $\tal_1$ is the constant from \cref{thm:ng:bounds}. On the other hand, arguing as in the proof of \cref{thm:converge}, we deduce that there exists $\tau_0'>0$ such that 
 \begin{align}
        2\nu\ZZ_N(t)&\leq 2\nu^3\left(\frac{\nu}{\bar{\mu}_1}\right)\left(\frac{\De\bar{\nu}_0}{\nu}\right)\de_0K_2^2,\label{est:ZN:0:second}\\
     |\lb B(w(t),w_N(t),Q_Nw(t))\rb|&\leq C\nu^3\left(\frac{\nu}{\bar{\mu}_1}\right)^{9/4}\left(\frac{\De\bar{\nu}_0}{\nu}\right)\de_0^2K_2^3,\label{est:B:0:second}\\
     |\lb (DB_Nu(t))w(t),w_N(t)\rb|
     %&\leq |(DB_Nu^{(1)})w^{(1)}||w_N^{(1)}|\leq C\nu(\s_1^{1/2}+G)^{3/2}\PP^{(1)}\notag\\
        &\leq C\nu^3(\de_0+1)^{3/2}(\s_1^{1/2}+\bar{G}_0)^{3/2}\left(\frac{\nu}{\bar{\mu}_1}\right)\left(\frac{\De\nu_0}{\nu}\right)\de_0K_2^2,\label{est:DB:0:second}
    \end{align}
for all $t\geq \tau_0'$, where we additionally applied \eqref{cond:mu:N:nu}. Recall that \eqref{est:EN:0:second}--\eqref{est:B:0:second} hold provided that $\bar{\mu}_1, N_1$ satisfy \eqref{cond:mu:N:tnu}, \eqref{cond:mu:N:nu}, \eqref{cond:mu:ng:H2}, $\eqref{cond:mu:ng:Hk}_3$, where $\tnu=\bar{\nu}_0$, $\mu=\bar{\mu}_1$. Hence, let us choose $t_1>t_1'\geq \tau_0\vee\tau_0'$, so that \eqref{est:EN:0:second}--\eqref{est:DB:0:second} all hold over the interval $J_1$.

Now choose $\be\in(0,1\wedge\de_0^{-1})$. Then we choose $\bar{\mu}_1\geq\bar{\nu}_0$ such that
    \begin{align}\label{cond:mu:proof:second}
    \frac{C}{\bar{\veps}_1}\frac{\de_0}{|1-\de_0|}\left(\frac{\bar{\nu}_0}{\bar{\mu}_1}\right)\left(\tal_1^2\w_1(1+\de_0)^4\bar{G}_0^2+\left(1+\frac{\de_0}{|1-\de_0|^{9/4}}K_2+\de_0(\s_1^{1/2}+\bar{G}_0)^{3/2}\right)K_2^2\right)\leq \be,
    \end{align}
where $C$ is sufficiently large (determined from the estimates above). Note that from \eqref{eq:nu:nu0:equiv}, \eqref{est:G:G0:equiv}, \eqref{est:tG0}, \eqref{eq:nu:nu0:equiv:second}, \eqref{est:G:G0:equiv:second}, \eqref{est:tG0:second}, $K_2$ can be replaced by a quantity that depends only $\al,\de_0,\s_1,\s_2,\bar{G}_0$. From our choice of $\be$, it follows that the analog of \eqref{eq:positivity} holds. Furthermore, from our choice of $\bar{\mu}_1$, it follows that
    \begin{align}\notag
        |\De\bar{\nu}_1|\leq \be|\De\bar{\nu}_0|.
    \end{align}
Making use of \eqref{eq:nu:nu0:equiv}, \eqref{est:G:G0:equiv}, \eqref{est:tG0}, \eqref{eq:nu:nu0:equiv:second}, \eqref{est:G:G0:equiv:second}, \eqref{est:tG0:second}, observe that the conditions \eqref{cond:mu:N:tnu}, \eqref{cond:mu:N:nu}, \eqref{cond:mu:ng:H2}, $\eqref{cond:mu:ng:Hk}_3$, \eqref{cond:mu:proof:second} may be summarized as
    \begin{align}\notag
        \bar{\gam}_1\leq\frac{\bar{\mu}_1}{\bar{\nu}_0}\leq\bar{\gam}N_1^2,
    \end{align}
for some constants $\bar{\gam}_1,\gam$, where $\gam$ depends only on $\al,\de_0,\s_1,\s_2,\bar{G}_0$ and $\bar{\gam}_1$ depends additionally on $\bar{\veps}_1^{-1},\w_1$. This establishes the base case, $M=1$.

Now suppose that for some $M\geq2$, there exist times $0=t_0<\dots<t_{M-1}$, intervals $J_m=[t_m',t_m)\subset[t_{m-1},t_{m})$, for $m=1,\dots, M-1$, and constants $\bar{\gam},\bar{\gam}_1,\dots,\bar{\gam}_{M-1}>0$, where $\bar{\gam}$ depends on $\al,\de_0,\s_1,\s_2,\bar{G}_0$, and $\bar{\gam}_m$ depends additionally on $\bar{\veps}_m^{-1}, \w_m$, where
    \begin{align}\label{def:omega:m}
        \w_m^{-1}:=\bar{\nu}_0|J_m|,
    \end{align}
such that if 
    \begin{align}\notag
        \bar{\veps}_1,\dots,\bar{\veps}_{M-1}>0,    
    \end{align}
where $\bar{\veps}_m$ is given by \eqref{def:bar:eps}, and $\bar{\mu}_1,\dots,\bar{\mu}_{M-1}>0$, $N_1,\dots, N_{M-1}\geq1$ satisfy
    \begin{align}\notag
        \bar{\gam}_m\leq\frac{\bar{\mu}_m}{\bar{\nu}_0}\leq\bar{\gam}N_m^2,
    \end{align}
for all $m=1,\dots, M-1$, then there exists, $\be\in(0,1\wedge\de_0^{-1})$ such that for $\bar{\nu}_1,\dots, \bar{\nu}_{M-1}$ given by \eqref{def:update:second}, corresponding to $J_1,\dots, J_{M-1}$, respectively, one has
    \begin{align}\notag
        |\De\bar{\nu}_m|\leq\be|\De\bar{\nu}_{m-1}|,
    \end{align}
for all $m=1,\dots, M-1$, and $\bar{\nu}_1,\dots,\bar{\nu}_{M-1}>0$. Note that the bar analogs of \eqref{eq:iterate} and \eqref{est:nu:num:equiv} also hold. Thus, in light of \cref{rmk:equiv} and \eqref{eq:nu:nu0:equiv}, \eqref{est:G:G0:equiv}, \eqref{est:tG0}, \eqref{eq:nu:nu0:equiv:second}, \eqref{est:G:G0:equiv:second}, \eqref{est:tG0:second}, we see that any dependence on $\nu$ appearing in \eqref{cond:mu:N:nu}, \eqref{cond:mu:ng:H2}, $\eqref{cond:mu:ng:Hk}_3$ may be replaced by $\bar{\nu}_0$ up to a pre-factor that depends only on $\be,\de_0$. In particular, all instances of $\tG,G$, where $\tG$ is defined by \eqref{def:tG} with $\tnu=\bar{\nu}_m$, can be replaced by $\bar{G}_0$ up to a pre-factor that depends on $\be,\de_0$. Lastly, in \eqref{cond:mu:N:tnu}, $\tnu=\bar{\nu}_m$ may be replaced by $\nu_0$ up to a pre-factor that depends only on $\be,\de_0$.

We re-initialize the feedback control system over $\til{I}_{M-1}=[t_{M-1},\infty)$ with $\tu(t_{M-1}^-)\in B_k(\al R_k)$, where $\al\geq1$ is the same scaling factor used for re-initializing in all previous intervals $\til{I}_0,\til{I}_1,\dots,\til{I}_{M-2}$, with the viscosities $\bar{\nu}_1,\bar{\nu}_2,\dots, \bar{\nu}_{M-1}$. From \eqref{def:update:second}, we have
   \begin{align}
     |\De\bar{\nu}_{M}&|\leq \frac{1}{\bar{\veps}_M\nu^2}\left(\frac{|\EE_N(t_{M})-\EE_N(t_{M}')|}{|J_{M}|}\right.\notag\\
        &\left.+\frac{1}{|J_{M}|}\int_{J_{M}}2\nu\ZZ_N(s)+|\lb B(w(s),w_N(s)),Q_Nw(s))\rb|+|\lb (DB_Nu(s))w(s),w_N(s)\rb| ds\right).\notag
    \end{align}
We may now argue as in the base case, applying \eqref{est:sensitivity:L2} to obtain a time $\tau_{M-1}>t_{M-1}$ such that \eqref{est:EN:0:second}, for all $t\geq \tau_{M-1}$, and a time $\tau_{M-1}'>t_{M-1}$ such that \eqref{est:ZN:0}--\eqref{est:DB:0:second} hold for all $t\geq \tau_{M-1}'$. Now we choose $J_M=[t_M',t_M)$ so that $t_M>t_M'\geq\tau_{M-1}\vee\tau_{M-1}'$. Hence, \eqref{est:EN:0:second}--\eqref{est:DB:0:second} hold over $J_M$. 

Next, we choose $\bar{\mu}_M\geq\bar{\nu}_0$ so that \eqref{cond:mu:N:tnu}, \eqref{cond:mu:N:nu}, \eqref{cond:mu:ng:H2}, $\eqref{cond:mu:ng:Hk}_3$ hold with $\tnu=\bar{\nu}_M$, $\mu=\bar{\mu}_M$, and  \eqref{cond:mu:proof:second}  holds with $\bar{\mu}_M$ replacing $\bar{\mu}_1$ and $\w_M$ replacing $\w_1$, where $\w_M$ is defined by \eqref{def:omega:m}. As in the base case, we invoke  \eqref{eq:nu:nu0:equiv}, \eqref{est:G:G0:equiv}, \eqref{est:tG0}, \eqref{eq:nu:nu0:equiv:second}, \eqref{est:G:G0:equiv:second}, \eqref{est:tG0:second}, as well as the bar analog of \eqref{est:nu:num:equiv}, in conjunction with \cref{rmk:equiv}, so that we may replace all dependencies on $\bar{\nu}_M, \nu, G$ by $\bar{\nu}_0,\de_0, \bar{G}_0$. %Moreover, the constant $C>0$ appearing in \eqref{cond:mu:proof:second}, may taken to be same as the one appearing in all of the previous instances, $m=1,\dots, M-1$.

Ultimately, due to this choice of $J_M$ and $\bar{\mu}_M$, we may apply the bounds \eqref{est:EN:0:second}--\eqref{est:DB:0:second} over $J_M$ to estimate $|\De\bar{\nu}_M|$ and arrive at
    \begin{align}\notag
        |\De\bar{\nu}_M|\leq \be|\De\bar{\nu}_{M-1}|.
    \end{align}
By iterating, we again argue that $\nu_M>0$, and notice that the conditions (with the notation adjusted accordingly) \eqref{cond:mu:N:tnu}, \eqref{cond:mu:N:nu}, \eqref{cond:mu:ng:H2}, $\eqref{cond:mu:ng:Hk}_3$,  and \eqref{cond:mu:proof:second} can be summarized as
    \begin{align}\notag
        \bar{\gam}_M\leq\frac{\bar{\mu}_M}{\bar{\nu}_0}\leq\bar{\gam}N_M^2,
    \end{align}
where $\bar{\gam}$ depends on $\al,\de_0,\s_1,\s_2,\bar{G}_0$ and  $\bar{\gam}_M$ depends additionally on $\bar{\veps}_M^{-1},\w_M$, as desired. This completes the proof.
\end{proof}

\section*{Acknowledgments}
The work of V.R.M. was partially supported by the PSC-CUNY
Research Award Program under grant PSC-CUNY 64335-00 52.  The author would like to thank the ADAPT group for their inspiration, encouragement, and comments throughout the course of this work. Lastly, the author would like to thank the referees for the insightful comments they shared to improve the manuscript.
\appendix

\section{Proof of Theorem \ref{thm:ng:bounds}}\label{sect:app:apriori}

\subsection{\texorpdfstring{$H^1$}{H1}--estimates}

We will first prove \eqref{def:rad:H1:ng}. Indeed, we take the scalar product of \eqref{eq:nse:ng:ff} with $A\tu$ in $L^2$, then invoke self-adjointness of $A$, to obtain
    \begin{align}\notag
        \frac{1}2\frac{d}{dt}\lVert\tu\rVert^2+\tnu|A\tu|^2+\mu\lVert \tu\rVert^2=\lb f,A\tu\rb+\mu\lVert Q_N\tu\rVert^2+\mu\lb A^{1/2}P_Nu,A^{1/2}P_N\tu\rb.
    \end{align}
The right-hand side may be estimated with the Cauchy-Schwarz inequality, Young's inequality, and the inverse Poincar\'e inequality, so that we may ultimately arrive at
    \begin{align}\label{est:enstrophy:ng}
        \frac{d}{dt}\lVert\tu\rVert^2+\tnu|A\tu|^2+\mu\lVert \tu\rVert^2\leq \frac{|f|^2}{\tnu}+\tnu\frac{\mu}{\tnu N^2}|A\tu|^2+\mu\lVert P_Nu\rVert^2.
    \end{align}
Given that $\mu, N, \tnu$ satisfy \eqref{cond:mu:N:tnu}, it then follows from \eqref{def:rad:H1} that
    \begin{align}
        \frac{d}{dt}\lVert\tu\rVert^2+\mu\lVert \tu\rVert^2\leq \nu^2\mu\left[\left(\frac{\nu}{\tnu}\right)\left(\frac{\nu}{\mu}\right) +2\right]G^2.\notag
    \end{align}
An application of Gronwall's inequality then yields
    \begin{align}\notag
        \lVert\tu(t)\rVert^2\leq e^{-\mu t}\lVert \tu_0\rVert^2+2\nu^2\tG^2(1-e^{-\mu t}),
    \end{align} 
for all $t\geq0$, where $\tG$ is defined by \eqref{def:tG}.

Furthermore, if $\tu_0\in B_1(\al R_1)$, for some $\al\geq1$, where $R_1$ is given by \eqref{def:rad:H1}, it follows that
    \begin{align}\notag
         \lVert\tu(t)\rVert^2\leq 2\nu^2\left[\al^2 +\left(\frac{\nu}{\tnu}\right)\left(\frac{\nu}{\mu}\right) +1\right]G^2,
    \end{align}
for all $t\geq0$. Upon setting 
    \[
        \tal_1^2:=4\al^2,
    \]
then recalling that $\al\geq1$, we obtain \eqref{def:rad:H1:ng}, as desired.
    
\subsection{\texorpdfstring{$H^2$}{H2}--estimates}
Upon taking the scalar product of \eqref{eq:nse:ng:ff} with $A^2\tu$ in $L^2$, then integrating by parts, we obtain
    \begin{align}\label{eq:palinstrophy:balance:ng}
        \frac{1}2\frac{d}{dt}|A\tu|^2&+\tnu|A^{3/2}\tu|^2+\mu|A\tu|^2\notag\\
        &=-\lb B(\tu,\tu),A^2\tu\rb+\lb f,A^2\tu\rb+\mu|AQ_N\tu|^2+\mu\lb AP_Nu, AP_N\tu\rb.
    \end{align}
We will first treat the trilinear term. Upon integrating by parts, we may rewrite it as
    \begin{align}
        \lb B(\tu,\tu),A^2\tu\rb&=\lb \tu^j\bdy_j\tu^k,\bdy_\ell^2\bdy_l^2\tu^k\rb=-\lb \bdy_\ell\tu^j\bdy_j\tu^k,\bdy_\ell\bdy_l^2\tu^k\rb-\lb\tu^j\bdy_j\bdy_\ell\tu^k,\bdy_\ell\bdy_l^2\tu^k\rb\notag\\
        &=\lb\bdy_\ell\bdy_l\tu^j\bdy_j\tu^k,\bdy_\ell\bdy_l\tu^k\rb+\lb\bdy_\ell\tu^j\bdy_j\bdy_l\tu^k,\bdy_\ell\bdy_l\tu^k\rb+\lb\bdy_\ell\tu^j\bdy_j\bdy_\ell\tu^k,\bdy_l^2\tu^k\rb.\notag
    \end{align}
It follows from H\"older's inequality, interpolation, the Poincar\'e inequality, \eqref{def:rad:H1:ng}, and Young's inequality that
    \begin{align}\label{est:BuuA2u}
        |\lb B(\tu,\tu),A^2\tu\rb|&\leq C\lVert\tu\rVert\left(\sum_{\ell,l}\Sob{\bdy_\ell\bdy_l\tu}{L^4}^2+\sum_{j,\ell,l}\Sob{\bdy_j\bdy_l\tu}{L^4}\Sob{\bdy_\ell\bdy_l\tu}{L^4}+\sum_{j,\ell}\Sob{\bdy_j\bdy_\ell\tu}{L^4}\Sob{A\tu}{L^4}\right)\notag\\
        &\leq C\lVert\tu\rVert|A^{3/2}\tu||A\tu|\leq C\nu\tal_1^2\left(\frac{\nu}{\tnu}\right)\tG^2|A\tu|^2+\frac{\tnu}{100}|A^{3/2}\tu|^2.
    \end{align}
We estimate the last three terms of \eqref{eq:palinstrophy:balance:ng} similar to \eqref{est:enstrophy:ng} and apply the result to \eqref{eq:palinstrophy:balance:ng} in combination with \eqref{est:BuuA2u} to arrive at
    \begin{align}\label{est:palinstrophy:ng}
    \frac{d}{dt}|A\tu|^2&+\tnu|A^{3/2}\tu|^2+\frac{3\mu}2|A\tu|^2\\
    &\leq C\nu\tal_1^2\left(\frac{\nu}{\tnu}\right)\tG^2|A\tu|^2+\s_1^2\nu^3\left(\frac{\nu}{\tnu}\right)G^2+\frac{\mu}{N^2\tnu}\tnu|A^{3/2}\tu|^2+C\mu\nu^2\left(\s_1^{1/2}+G\right)^2G^2,\notag
    \end{align}
where we applied \eqref{def:rad:Hk} with $k=2$ in estimating the last term in \eqref{eq:palinstrophy:balance:ng}. Since $\mu, N$ satisfy \eqref{cond:mu:N:tnu} and \eqref{cond:mu:ng:H2}, for $k=2$, where $\al_2^2=100C\tal_1^2$, we may reduce \eqref{est:palinstrophy:ng} to 
    \begin{align}\notag
        \frac{d}{dt}|A\tu|^2+\mu|A\tu|^2\leq  C\mu\nu^2\s_1^2\left(\s_1^{1/2}+G\right)^2\tG^2,
    \end{align}
An application of Gronwall's inequality then yields
    \begin{align}\notag
        |A\tu(t)|^2\leq e^{-\mu t}|A\tu_0|^2+C\nu^2\s_1^2\left(\s_1^{1/2}+G\right)^2\tG^2.
    \end{align}
Thus, assuming that $\tu_0\in B_2(\al R_2)$, for some $\al\geq1$, where $R_2$ is given by \eqref{def:rad:Hk}, where $k=2$, it follows that
    \begin{align}\notag
        |A\tu(t)|^2\leq C(\al^2+\s_1^2)\nu^2\left(\s_1^{1/2}+G\right)^2\tG^2,
    \end{align}
for all $t\geq0$. Note that we used the fact that $G\leq\tG$. Upon setting     
    \[
        \tal_2^2:=C(\al^2+\s_1^2),
    \]
we obtain \eqref{def:rad:H2:ng}, as desired.

\subsection{\texorpdfstring{$H^k$}{Hk}--estimates}
First, observe that by Placherel's theorem, the Cauchy-Schwarz inequality, and the Poincar\'e inequality one sees that $|A^{k/2}\phi|^2\leq C(|\bdy_1^k\tu|^2+|\bdy_2^k\tu|^2)$. Thus, it suffices to obtain estimates on $\bdy_\ell^k\tu$ for all $k\geq3$ and $\ell=1,2$. In this regard, we will make use of the following notation:
    \[
        \Sob{\phi}{\dot{H}^k}^2=\sum_{\ell=1,2}|\bdy_\ell^k\phi|^2.
    \]

Upon taking the scalar product of \eqref{eq:nse:ng:ff}, with $(-1)^k\bdy_\ell^{2k}$, we obtain 
    \begin{align}\label{eq:balance:Hk:ng}
        &\frac{1}2\frac{d}{dt}|\bdy_\ell^k\tu|^2+\tnu|\nabla\bdy_\ell\tu|^2+\mu|\bdy_\ell^k\tu|^2\notag\\
        &=-(-1)^k\lb B(\tu,\tu),\bdy_\ell^{2k}\tu\rb+(-1)^k\lb f,\bdy_\ell^{2k}\tu\rb+\mu|\bdy_\ell^{k}Q_N\tu|^2+\mu\lb \bdy_\ell^kP_Nu,\bdy_\ell^kP_N\tu\rb.
    \end{align}
First, observe that the trilinear term may be expanded using integration by parts, the Leibniz rule, and \eqref{eq:ortho:B}, to obtain
    \begin{align}\label{eq:Buu:Hku}
    &-(-1)^k\lb B(\tu,\tu),\bdy_\ell^{2k}\tu\rb=-\lb \bdy_\ell^kB(\tu,\tu),\bdy_\ell^k\tu\rb\notag\\
    &=-\left(k\lb\bdy_\ell\tu^j\bdy_j\bdy_\ell^{k-1}\tu^i,\bdy_\ell^k\tu^i\rb+\lb\bdy_\ell^k\tu^j\bdy_j\tu^i,\bdy_\ell^k\tu^i\rb+\sum_{l=2}^{k-1}c_{l,k}\lb\bdy_\ell^l\tu^j\bdy_j\bdy_\ell^{k-l}\tu^i,\bdy_\ell^k\tu^i\rb\right).
    \end{align}
By H\"older's inequality and interpolation we have
    \begin{align}
        |\lb\bdy_\ell\tu^j\bdy_j\bdy_\ell^{k-1}\tu^i,\bdy_\ell^k\tu^i\rb|&\leq \lVert\tu\rVert\Sob{\bdy_j\bdy_\ell^{k-1}\tu}{L^4}\Sob{\bdy_\ell^k\tu}{L^4}\leq C\lVert\tu\rVert|A^{(k+1)/2}\tu||A^{k/2}\tu|\notag\\
         |\lb\bdy_\ell^k\tu^j\bdy_j\tu^i,\bdy_\ell^k\tu^i\rb|&\leq C\lVert\tu\rVert|A^{(k+1)/2}\tu||A^{k/2}\tu|.\notag
    \end{align}
Similarly, for $2\leq l\leq k-1$, we have
    \begin{align}
       |\lb\bdy_\ell^l\tu^j\bdy_j\bdy_\ell^{k-l}\tu^i,\bdy_\ell^k\tu^i\rb|&\leq \Sob{\bdy_\ell^l\tu}{L^\infty}|\bdy_j\bdy_\ell^{k-l}\tu||A^{k/2}\tu|\notag\\
       &\leq C|A^{(l+2)/2}\tu|^{1/2}|A^{l/2}\tu|^{1/2}|A^{(k-l+1)/2}\tu||A^{k/2}\tu|\notag\\
       &\leq C|A^{(k+1)/2}\tu|^{\frac{k-2}{k-1}}|A\tu|^{\frac{k}{k-1}}|A^{k/2}\tu|\notag.
    \end{align}
Returning to \eqref{eq:Buu:Hku}, then applying Young's inequality, Poincar\'e's inequality, \eqref{def:rad:H1:ng}, and \eqref{def:rad:H2:ng} with $k=2$, we obtain
    \begin{align}
       &\sum_\ell |\lb B(\tu,\tu),\bdy_\ell^{2k}\tu\rb|\leq C\lVert\tu\rVert|A^{(k+1)/2}\tu||A^{k/2}\tu|+C|A^{(k+1)/2}\tu|^{\frac{k-2}{k-1}}|A\tu|^{\frac{k}{k-1}}|A^{k/2}\tu|\notag\\
       &\leq\frac{\tnu}{100}|A^{(k+1)/2}\tu|^2+C\nu\left(\frac{\nu}{\tnu}\right)\left(\frac{\lVert\tu\rVert}{\nu}\right)^2|A^{k/2}\tu|^2+C\frac{|A\tu|^2}{\tnu^{\frac{k-2}{k}}}|A^{k/2}\tu|^{\frac{2(k-1)}{k}}\notag\\
       &\leq \frac{\tnu}{100}|A^{(k+1)/2}\tu|^2+C\nu\left(\frac{\nu}{\tnu}\right)\left(\frac{\lVert\tu\rVert}{\nu}\right)^2|A^{k/2}\tu|^2
       \notag\\
       &\quad+C\nu^2\mu\left(\frac{\nu}{\tnu}\right)^{k-2}\left(\frac{\nu}{\mu}\right)^k\left(\frac{|A\tu|}{\nu}\right)^{2k}+\frac{\mu}{100}|A^{k/2}\tu|^2\label{est:Buu:Aku}\\
       &\leq \frac{\tnu}{100}|A^{(k+1)/2}\tu|^2+\left[C\tal_1^2\left(\frac{\nu}{\tnu}\right)\left(\frac{\nu}{\mu}\right)\tG^2+\frac{1}{100}\right]\mu|A^{k/2}\tu|^2\notag\\
       &\quad+C\nu^2\mu\left(\frac{\nu}{\tnu}\right)^{k-2}\left(\frac{\nu}{\mu}\right)^k\tal_2^2\left(\s_1^{1/2}+G\right)^k\tG^{k}\notag\\
       &\leq \frac{\tnu}{100}|A^{(k+1)/2}\tu|^2+\frac{\mu}{50}|A^{k/2}\tu|^2+C\nu^2\mu\left(\frac{\nu}{\tnu}\right)^{k-2}\left(\frac{\nu}{\mu}\right)^k\tal_2^2\left(\s_1^{1/2}+G\right)^k\tG^{k}.\notag
    \end{align}

Now we return to \eqref{eq:balance:Hk:ng} and sum over $\ell=1,2$. Similar to \eqref{est:enstrophy:ng}, we estimate the last three terms of \eqref{eq:balance:Hk:ng}, but additionally invoke \eqref{def:rad:Hk}, to obtain
    \begin{align}\label{est:Hk:ng}
    \frac{d}{dt}\Sob{\tu}{\dot{H}^k}^2&+\frac{\tnu}2\Sob{\nabla\tu}{\dot{H}^k}^2+\frac{3\mu}2\Sob{\tu}{\dot{H}^k}^2\notag\\
    &\leq  \frac{\Sob{f}{\dot{H}^{k-1}}^2}{\tnu}+\frac{\mu}{N^2\nu}\nu\Sob{\nabla\tu}{\dot{H}^k}^2+2\mu|A^{k/2}u|^2+C\nu^2\mu\left(\frac{\nu}{\tnu}\right)^{k-2}\left(\frac{\nu}{\mu}\right)^k\tal_2^2\left(\s_1^{1/2}+G\right)^k\tG^{k}\notag\\
    &\leq \nu^3\left(\frac{\nu}{\tnu}\right)\s_{k-1}^2G^2+\frac{\mu}{N^2\tnu}\tnu\Sob{\nabla\tu}{\dot{H}^k}^2\notag\\
    &\quad+C\nu^2\mu\left[(\s_{k-1}^{1/k}+G)^{2(k-1)}G^2+\tal_2^2\left(\frac{\nu}{\tnu}\right)^{k-2}\left(\frac{\nu}{\mu}\right)^k\left(\s_1^{1/2}+G\right)^k\tG^{k}\right],
    \end{align}
Then by \eqref{cond:mu:N:tnu}, \eqref{cond:mu:ng:H2}, \eqref{cond:mu:ng:Hk}, and Gronwall's inequality, it follows that
    \begin{align}
        \Sob{\tu(t)}{\dot{H}^k}^2&\leq e^{-\mu t}\Sob{\tu_0}{\dot{H}^k}^2+C\nu^2\left[\left(\frac{\nu}{\tnu}\right)\left(\frac{\nu}{\mu}\right)\s_{k-1}^2+\left(\s_{k-1}^{1/k}+G\right)^{2(k-1)}\right]G^2\notag\\
        &\leq e^{-\mu t}\Sob{\tu_0}{\dot{H}^k}^2+C\nu^2\s_{k-1}^{\frac{2}k}(\s_{k-1}^{1/k}+G)^{2(k-1)}\tG^2.\notag
    \end{align}
 Assuming that $\tu_0\in B_k(\al R_k)$, for some $\al\geq1$, where $R_k$ is given by \eqref{def:rad:Hk}, it follows that
    \begin{align}\notag
        |A^{k/2}\tu(t)|^2\leq C\nu^2(\al^2+\s_{k-1}^{2/k})(\s_{k-1}^{1/k}+G)^{2(k-1)}\tG^2,
    \end{align}
for all $t\geq0$ and $k\geq3$. Note that we invoked the facts that $G\leq\tG$ and $\s_{k-1}\geq1$. Upon setting
    \[
        \tal_k^2=C(\al^2+\s_{k-1}^{2/k}),
    \]
we obtain \eqref{def:rad:Hk:ng}, for $k\geq3$, as desired.

\section{Higher-order sensitivity-type bounds}\label{sect:app:sensitivity}

Suppose that $\nu,\tnu>0$ and consider the unique solutions $u,\tu$ of \eqref{eq:nse} and \eqref{eq:nse:ng} corresponding to $\nu,\tnu$, respectively. Let $w=\tu-u$. Then we may rewrite \eqref{eq:nse:diff:ng} as
    \begin{align}\label{eq:nse:diff:ng:equiv}
        \frac{dw}{dt}+\nu Aw+B(w,w)+(DBu)w=-(\De\nu)A\tu-\mu w+\mu Q_Nw,\quad \De\nu=\tnu-\nu.
    \end{align}
We proceed in a bootstrap fashion. In particular, we will first derive refined sensitivity bounds in $L^2$ and $H^1$, before proceeding to $H^2$ and higher.

\subsection{\texorpdfstring{$L^2$}{L2}--estimates}
Upon taking the $L^2$ scalar product of \eqref{eq:nse:diff:ng:equiv} with $w$, we obtain
    \begin{align}\label{eq:diff:balance:E}
        \EEd+2\nu\ZZ+2\mu\EE=-\lb (DBu)w,w\rb-(\De\nu)\lb A^{1/2}\tu,A^{1/2}w\rb+\mu|Q_Nw|^2=E_1+E_2+E_3,
    \end{align}
where we used the fact that $A^p$ is self-adjoint and the notation $\EE,\EEd,\ZZ$ introduced in \eqref{def:diff:energy}.

Observe that
    \begin{align}\notag
        E_1=-\lb B(w,u),w\rb.
    \end{align}
Then by H\"older's inequality, interpolation, Young's inequality, and \eqref{def:rad:H1} we have
    \begin{align}
        |E_1|\leq \Sob{w}{L^4}^2\lVert u\rVert\leq \lVert w\rVert|w|\lVert u\rVert\leq C\nu\left(\frac{\nu}{\mu}\right)\left(\frac{\lVert u\rVert}{\nu}\right)^2\ZZ+\frac{\mu}{100}\EE\leq C\nu\left(\frac{\nu}\mu\right) G^2\ZZ+\frac{\mu}{100}\EE\notag.
    \end{align}
For $E_2$, we estimate with the Cauchy-Schwarz inequality, Young's inequality, and \eqref{def:rad:H1:ng} to obtain
    \begin{align}\notag
        |E_2|\leq |\De\nu\lVert\tu\rVert\lVert w\rVert\leq C\left(\frac{(\De\nu)^2}{\nu}\right)\lVert\tu\rVert^2+\frac{\nu}{100}\EE\leq C\tal_1^2\nu^3\left(\frac{|\De\nu|}{\nu}\right)^2\tG^2+\frac{\nu}{100}\ZZ.
    \end{align}
For $E_3$, we estimate with the inverse Poincar\'e inequality and \eqref{cond:mu:N:nu} to obtain
    \begin{align}\notag
        E_3\leq \frac{\mu}{N^2\nu}\nu\lVert w\rVert^2\leq \frac{\nu}{2}\ZZ.
    \end{align}
    
Combining the estimates for $E_1$--$E_3$ and using the fact that $\mu$ satisfies \eqref{cond:mu:sensitivity1}, it follows that
    \begin{align}
        \EEd+\nu\ZZ+\mu\EE\leq C\tal_1^2\nu^3\left(\frac{|\De\nu|}{\nu}\right)^2\tG^2.\notag
    \end{align}
By Gronwall's inequality, we arrive at
    \begin{align}\label{est:sensitivity:L2}
        \EE(t)\leq e^{-\mu (t-t_0)}\EE(t_0)+C\tal_1^2\nu^2\left(\frac{\nu}{\mu}\right)\left(\frac{|\De\nu|}{\nu}\right)^2\tG^2\left(1-e^{-\mu(t-t_0)}\right),
    \end{align}
for any $t_0\geq0$. In particular
    \begin{align}\notag
         \EE(t)\leq e^{-\mu (t-t_0)}\EE(t_0)+\nu^2\left(\frac{\nu}{\mu}\right)\left(\frac{|\De\nu|}{\nu}\right)^2K_0^2,
    \end{align}
where
    \begin{align}\label{def:K0}
      K_0^2=C\tal_1^2\tG^2.
    \end{align}
    
\subsection{\texorpdfstring{$H^1$}{H1}--estimates}

Taking the $L^2$ scalar product of \eqref{eq:nse:diff:ng:equiv} with $Aw$, we obtain
    \begin{align}\notag
        \ZZd+2\nu\PP+2\mu\ZZ=-\lb (DBu)w,Aw\rb-(\De\nu)\lb A\tu, Aw\rb+\mu\lVert Q_Nw\rVert^2=Z_1+Z_2+Z_3.
    \end{align}
Observe that by integrating by parts, we may rewrite $Z_1$ as
    \begin{align}\notag
        Z_1=-\lb B(u,w),Aw\rb-\lb B(w,u),Aw\rb=\lb \bdy_ku^j\bdy_jw^\ell\bdy_kw^\ell\rb-\lb B(w,u),Aw\rb,
    \end{align}
where we adopt the convention of summing over repeated indices. We then estimate $K_1$ with H\"older's inequality, interpolation, \eqref{def:rad:H1}, and \eqref{def:rad:Hk} for $k=2$, and Poincare's inequality, we obtain
    \begin{align}
        |Z_1|&\leq \lVert u\rVert\Sob{\nabla w}{L^4}^2+\Sob{w}{L^4}\Sob{\nabla u}{L^4}|Aw|\notag\\
        &\leq   C\left(\lVert u\rVert\lVert w\rVert+\lVert w\rVert^{1/2}|w|^{1/2}|Au|^{1/2}\lVert u\rVert^{1/2}\right)|Aw|\notag\\
        &\leq C\frac{\lVert u\rVert^2}{\nu}\ZZ+C\frac{|Au|\lVert u\rVert}{\nu}\lVert w\rVert|w|+\frac{\nu}{100}|Aw|^2\notag\\
        &\leq C\nu G^2\ZZ+ C\nu{\left(\s_1^{1/2}+G\right)G}\ZZ^{1/2}\EE^{1/2}+\frac{\nu}{100}\PP\notag\\
        &\leq C\nu G^2\ZZ+C\nu\left(\s_1^{1/2}+G\right)^2\EE+\frac{\nu}{100}\PP\notag\\
        &\leq C\nu\left(\s_1^{1/2}+G\right)^2\ZZ+\frac{\nu}{100}\PP\notag.
    \end{align}
Next, we estimate $K_2$ with the Cauchy-Schwarz, Young's inequality, and \eqref{def:rad:Hk:ng} for $k=2$, to obtain
    \begin{align}
        |Z_2|&\leq |\De\nu||A\tu||Aw|\leq C\nu\left(\frac{|\De\nu|}{\nu}\right)^2|A\tu|^2+\frac{\nu}{100}\PP\notag\\
        &\leq C\nu^2\left(\frac{|\De\nu|}{\nu}\right)^2\tal_2^2\s_1\left(\s_1^{1/2}+G\right)^2\tG^2+\frac{\nu}{100}\PP.\notag
    \end{align}
Lastly for $K_3$, we estimate with the inverse Poincar\'e inequality and \eqref{cond:mu:N:nu} to obtain
    \begin{align}\notag
        Z_3\leq \frac{\mu}{N^2\nu}\nu|Aw|^2\leq \frac{\nu}2\PP.
    \end{align}

Combining the estimates for $Z_1$--$Z_3$ and invoking \eqref{cond:mu:sensitivity1} holds, it follows that
    \begin{align}\notag
        \ZZd+\nu\PP+\mu\ZZ\leq C\tal_2^2\nu^2\left(\frac{|\De\nu|}{\nu}\right)^2\s_1\left(\s_1^{1/2}+G\right)^2\tG^2.
    \end{align}
Then by Gronwall's inequality, we have
    \begin{align}\label{est:sensitivity:H1}
        \ZZ(t)\leq e^{-\mu(t-t_0)}\ZZ(t_0)+C\tal_2^2\nu^2\left(\frac{\nu}{\mu}\right)\left(\frac{|\De\nu|}{\nu}\right)^2\left(\s_1^{1/2}+\tG\right)^2\tG^2\left(1-e^{-\mu(t-t_0)}\right)
    \end{align}
for all $t\geq t_0\geq 0$. In particular
    \begin{align}\notag
        \ZZ(t)\leq e^{-\mu(t-t_0)}\ZZ(t_0)+\nu^2\left(\frac{\nu}{\mu}\right)\left(\frac{|\De\nu|}{\nu}\right)^2K_1^2,
    \end{align}
where
    \begin{align}\label{def:K1}
        K_1^2:=C\tal_2^2\left(\s_1^{1/2}+\tG\right)^2\tG^2.
    \end{align}

\subsection{\texorpdfstring{$H^2$}{H2}--estimates}
Taking the $L^2$ scalar product of \eqref{eq:nse:diff:ng:equiv} with $A^2w$, we obtain
    \begin{align}
        \PPd+\nu|A^{3/2}w|^2+2\mu\PP&=-\lb B(w,w),A^2w\rb-\lb (DBu)w,A^2w\rb-\De\nu\lb A\tu, A^2w\rb+\mu|Q_NAw|^2\notag\\
        &=P_0+P_1+P_2+P_3.\notag
    \end{align}
We estimate $P_0$ as in \eqref{est:BuuA2u} and obtain
    \begin{align}\notag
        |P_0|\leq C\lVert w\rVert|A^{3/2}w||Aw|\leq \frac{C}{\nu}\ZZ\PP+\frac{\nu}{100}|A^{3/2}w|^2.
    \end{align}
For $P_1$, we integrate by parts first to write
    \begin{align}
       & \lb B(u,w),A^2w\rb=\lb u^j\bdy_jw^k,\bdy_\ell^2\bdy_l^2w^k\rb=-\lb \bdy_\ell u^j\bdy_jw^k,\bdy_\ell\bdy_l^2w^k\rb\notag\\
        &=\lb\bdy_\ell\bdy_l u^j\bdy_jw^k,\bdy_\ell\bdy_l w^k\rb+\lb\bdy_\ell u^j\bdy_j\bdy_lw^k,\bdy_\ell\bdy_lw^k\rb=\lb B(\bdy_\ell\bdy_lu,w)+B(\bdy_\ell u,\bdy_lw),\bdy_\ell\bdy_lw\rb,\notag
    \end{align}
and similarly
    \begin{align}
       & \lb B(w,u),A^2w\rb=\lb w^j\bdy_j w^k,\bdy_\ell^2\bdy_l^2w^k\rb=-\lb \bdy_\ell w^j\bdy_j u^k,\bdy_\ell\bdy_l^2w^k\rb-\lb w^j\bdy_j\bdy_\ell u^k,\bdy_\ell\bdy_l^2w^k\rb\notag\\
       % &=\lb\bdy_\ell\bdy_l w^j\bdy_ju^k,\bdy_\ell\bdy_lw^k\rb+\lb\bdy_\ell w^j\bdy_j\bdy_lu^k,\bdy_\ell\bdy_lw^k\rb+\lb\bdy_lw^j\bdy_j\bdy_\ell u^k,\bdy_\ell\bdy_lw^k\rb+\lb w^j\bdy_j\bdy_\ell\bdy_lu^k,\bdy_\ell\bdy_lw^k\rb\notag\\
        &=\lb B(\bdy_\ell\bdy_lw,u)+ B(\bdy_\ell w,\bdy_lu)+B(\bdy_lw,\bdy_\ell u)+ B(w,\bdy_\ell\bdy_l u),\bdy_\ell\bdy_lw\rb.\notag
    \end{align}
Thus by H\"older's inequality, interpolation, and Young's inequality, we have
    \begin{align}
        |\lb B(u,w),A^2w\rb|&\leq \Sob{\bdy_\ell\bdy_lu}{L^4}\Sob{\nabla w}{L^4}|Aw|+\lVert u\rVert\Sob{\bdy_j\bdy_lw}{L^4}\Sob{\bdy_\ell\bdy_lw}{L^4}\notag\\
        &\leq C|A^{3/2}u|^{1/2}|Au|^{1/2}|Aw|^{3/2}\lVert w\rVert^{1/2}+C\lVert u\rVert|A^{3/2}w||Aw|\notag\\
        &\leq C\nu\left(\frac{\nu}{\mu}\right)^3\left(\frac{|A^{3/2}u||Au|}{\nu^2}\right)^2\ZZ+\frac{\mu}{100}\PP+C\nu\left(\frac{\lVert u\rVert}{\nu}\right)^2\PP+\frac{\nu}{100}|A^{3/2}w|^2.\notag
    \end{align}
Similarly
    \begin{align}
    &|\lb B(w,u),A^2w\rb|\notag\\
    &\leq \lVert u\rVert\Sob{\bdy_\ell\bdy_lw}{L^4}\Sob{\bdy_\ell\bdy_lw}{L^4}+\Sob{\bdy_j\bdy_lu}{L^4}\Sob{\nabla w}{L^4}\Sob{\bdy_\ell\bdy_lw}{L^4}+\Sob{w}{L^\infty}|A^{3/2}u||Aw|\notag\\
    &\leq C\lVert u\rVert|A^{3/2}w||Aw|+C|A^{3/2}u|^{1/2}|Au|^{1/2}|Aw|^{3/2}\lVert w\rVert^{1/2}+C|w|^{1/2}|A^{3/2}u||Aw|^{3/2}\notag\\
    &\leq C\nu\left(\frac{\nu}{\mu}\right)^3\left(\frac{|A^{3/2}u||Au|}{\nu^2}\right)^2\ZZ+\frac{\mu}{100}\PP+C\nu\left(\frac{\lVert u\rVert}{\nu}\right)^2\PP+C\nu\left(\frac{\nu}{\mu}\right)^3\left(\frac{|A^{3/2}u|}{\nu}\right)^4\EE\notag\\
    &\quad+\frac{\nu}{100}|A^{3/2}w|^2\notag.
    \end{align}
Hence, upon combining these estimates, it follows from \eqref{def:rad:Hk} and Poincare's inequality that
    \begin{align}
    |P_1|\leq& C\mu\left(\frac{\nu}{\mu}\right)^4(\s_2^{1/3}+G)^4(\s_1^{1/2}+G)^2G^4\ZZ+C\nu G^2\PP+C\mu\left(\frac{\nu}{\mu}\right)^4(\s_2^{1/3}+G)^8G^4\EE\notag\\
    &\quad+\frac{\mu}{50}\PP+\frac{\nu}{50}|A^{3/2}w|^2.\notag
    \end{align}
Finally, for $P_2, P_3$ we apply the Cauchy-Schwarz inequality, Young's inequality, the inverse Poinca\'re inequality, \eqref{def:rad:Hk:ng}, and \eqref{cond:mu:N:nu} to estimate
    \begin{align}\notag
        |P_2|+|P_3|&\leq |\De\nu||A^{3/2}\tu||A^{3/2}w|+\frac{\mu}{N^2\nu}\nu|A^{3/2}w|^2\notag\\
        &\leq C\nu\left(\frac{|\De\nu|}{\nu}\right)^2|A^{3/2}\tu|^2+\nu\left(\frac{\mu}{N^2\nu}+\frac{1}{100}\right)|A^{3/2}w|^2\notag\\
        &\leq C\nu^3\left(\frac{|\De\nu|}{\nu}\right)^2\tal_2^2\s_2^{2/3}(\s_2^{1/3}+G)^4\tG^2+\frac{\nu}{50}|A^{3/2}w|^2.\notag
    \end{align}

Combining the estimates for $P_0$--$P_3$, we obtain
    \begin{align}\notag
        \frac{d}{dt}\PP+\frac{3\mu}2\PP&\leq \frac{C}{\nu}\ZZ\PP+C\mu\left(\frac{\nu}{\mu}\right)^4(\s_2^{1/3}+G)^4G^4\left[(\s_1^{1/2}+G)^2\ZZ+(\s_2^{1/3}+G)^4\EE\right]\notag\\
        &\quad+C\nu^3\tal_2^2\s_2^{2/3}(\s_2^{1/3}+G)^4\tG^2\left(\frac{|\De\nu|}{\nu}\right)^2.\notag
    \end{align}
Now observe that by the Poincar\'e inequality, we have $\EE,\ZZ\leq\PP$. Hence, upon invoking \eqref{cond:mu:sensitivity1}, we arrive at
    \begin{align}
        &\frac{d}{dt}\PP+\mu\PP\leq \frac{C}{\nu}\ZZ\PP+C\nu^3\tal_2^2\s_2^{2/3}(\s_2^{1/3}+G)^4\tG^2\left(\frac{|\De\nu|}{\nu}\right)^2.\notag
    \end{align}
Now choose $t_0'>t_0$, sufficiently large so that
    \begin{align}\label{eq:PP:large:time}
        e^{-\mu(t_0'-t_0)}\ZZ(t_0)\leq C\nu\left[(\s_1^{1/2}+G)^2+(\s_2^{1/3}+G)^4\right]^{1/4}(\s_2^{1/3}+G)G,
    \end{align}
By Gronwall's inequality, \eqref{est:sensitivity:H1}, \eqref{cond:mu:sensitivity1}, \eqref{cond:mu:sensitivity2}, it now follows that
    \begin{align}\label{est:sensitivity:H2}
        \PP(t)\leq& e^{-\mu(t-t_0')}\PP(t_0') +\nu^2\left(\frac{\nu}{\mu}\right)\left(\frac{|\De\nu|}{\nu}\right)^2K_2,
    \end{align}
for all $t\geq t_0'$, where 
    \begin{align}\notag
        K_2^2:=2C\tal_2^2\s_2^{2/3}(\s_2^{1/3}+G)^4\tG^2.
    \end{align}

%\begin{footnotesize}

\providecommand{\href}[2]{#2}

%\end{footnotesize}

%\begin{multicols}{1}

\vfill 

\noindent Vincent R. Martinez\\
{\footnotesize
Department of Mathematics \& Statistics\\
CUNY Hunter College\\
Web: \url{http://math.hunter.cuny.edu/vmartine/}\\
Email: \url{vrmartinez@hunter.cuny.edu}\\
}

%\end{multicols}

\end{document}